\documentclass[11pt,leqno]{amsart}
\usepackage{a4wide}
\usepackage{arydshln}
\usepackage[width=17cm]{geometry}
\usepackage{placeins}
\usepackage{hyperref}
\hypersetup{pdftex,colorlinks=true,allcolors=blue}
\usepackage{hypcap}
\usepackage{pgf,tikz}
\usepackage[onelanguage,vlined]{algorithm2e}
\usepackage{listings}
\usepackage[onelanguage,vlined]{algorithm2e}
\usepackage{amssymb}
\usepackage{amsmath}
\usepackage{gensymb}% angle degree symbol

\newcommand{\CC}{\mathbb{C}}
\newcommand{\FF}{\mathbb{F}}

\newcommand{\QQ}{\mathbb{Q}}

\newcommand{\calH}{\mathcal{H}}

\newcommand{\calB}{\mathcal{B}}
\newcommand{\calS}{\mathcal{S}}

\newcommand{\PGL}{\operatorname{PGL}}
\newcommand{\GL}{\operatorname{GL}}
\newcommand{\SO}{\operatorname{SO}}

\newcommand{\End}{\operatorname{End}}

\usepackage{amsthm}

\newtheorem*{theorem*}{Theorem}

\newtheorem*{question*}{Question}
\theoremstyle{definition}
\theoremstyle{definition}

\title[On refined Bruhat decompositions and endomorphism 
algebras of Gelfand-Graev rep.]{On refined Bruhat decompositions and endomorphism 
algebras of Gelfand-Graev representations}\author{Alessandro Paolini}
\address{
Technische Universit\"at Kaiserslautern\\
Fachbereich Mathematik\\
Postfach 3049\\
67653 Kaiserslautern\\
}
\email{paolini@mathematik.uni-kl.de}
\author{Iulian I. Simion}
\address{
Babe\c s-Bolyai University\\
Department of Mathematics\\
Mathematica 9\\
400157, Ploie\c sti 23-25, Cluj-Napoca\\
}
\email{simion@math.ubbcluj.ro}

\begin{document}

\begin{abstract}
Let $G$ be a finite reductive group defined over $\FF_q$, with $q$ a power of a prime $p$. 
Motivated by a problem recently posed by C. Curtis, 
we first develop an algorithm to 
express each element of $G$ into a canonical form in terms of 
a refinement of a Bruhat decomposition, 
and we then use the output of the algorithm to explicitly determine the structure constants of the endomorphism algebra 
of a %standard 
Gelfand-Graev representation of $G$ %in the cases 
when $G=\PGL_3(q)$ for an arbitrary prime $p$, and when $G=\SO_5(q)$ for $p$ odd. 
 
\end{abstract}

\thanks{\hspace*{-1.3em}Date: \today. 
The first author acknowledges financial support from the SFB-TRR 195.
}

\maketitle
\setcounter{tocdepth}{1}

Let $G$ be the fixed point subgroup $\bar G^{F}$ of a reductive algebraic group $\bar G$ under a Frobenius endomorphism $F$. A special role in the representation theory of finite reductive groups is played by Gelfand-Graev representations. These are certain representations induced from a linear character of a maximal unipotent subgroup of $G$. If the center of the ambient algebraic group $\bar G$ is connected, then a Gelfand-Graev representation $\Gamma$ is unique up to conjugacy and the irreducible components of $\Gamma$ were described in \cite{DL76}. In general, these representations were parametrized and decomposed in \cite[Section 14]{DM91} and \cite[Section 3]{DLM92}.

Denote by $\calH$ the Hecke algebra (that is, the $G$-endomorphism algebra) of the module affording $\Gamma$. 
Since $\Gamma$ is multiplicity free \cite[Theorem 8.1.3]{carter_finite}, the algebra $\calH$ is abelian. 
C. Curtis parametrized in \cite{curtis93} the irreducible representations of $\calH$ by pairs $(T,\theta)$ where 
$T=\bar T^{F}$ for an $F$-stable maximal torus $\bar T$ of $\bar G$ and where $\theta$ is an irreducible character of $T$. 
%In one of his most recent works \cite{curtis}, he poses the following remark
Each irreducible representation $f_{T,\theta}$ of the algebra $\calH$ is shown to have the following factorization \cite[Theorem 4.2]{curtis93}, 
$$
f_{T,\theta}=\hat\theta\circ f_T,
$$
where $f_T:\calH\rightarrow\bar\QQ_\ell T$ is a homomorphism of algebras and $\hat\theta$ is the linear extension of $\theta$ to the group algebra $\bar\QQ_\ell T$. The homomorphism $f_T$ is independent of $\theta$. Curtis' homomorphism $f_T$ was considered in \cite{BK08} in the context of $\ell$-modular Gelfand-Graev representations with $\ell \ne p$. %FIX! In the modular case, 
Here it is shown that if $\ell$ does not divide the order of the Weyl group of $G$, then the behavior of such endomorphism algebras is 
generic along all prime powers $q$. %similar for all primes 

In order to obtain a constructive description of each $f_T$, one needs to consider structure constants of $\calH$ with respect to some basis. 
Building on work of Kawanaka \cite{kawanaka1975unipotent} and Deodhar \cite{deodhar85}, Curtis described in \cite{curtis} an algorithm for obtaining structure constants for $\calH$ with respect to a standard basis %which is
parametrized by certain elements in $N_G(T)$. In general, the determination of the structure constants requires a deep understanding of the following intersections of shifted $U$-double cosets, 
$$
Ux U\cap zUy^{-1}U,
\qquad\text{for}\qquad U \in \mathrm{Syl}_p(G) \qquad \text{and} \qquad  x,y,z\in N_G(\bar T). 
$$
A weaker form of these intersections are the \emph{refined Bruhat cells} which are obtained by replacing $U$ with $B=N_G(U)$, i.e. $Bx B\cap zBy^{-1}B$.
Here Curtis also raises the problem 
of explicitly obtaining formulas for the structure constants, which can be used to give a combinatorial proof of the existence of the homomorphisms $f_T$, see \cite[Section 4]{curtis}. 

The interest in the above cell intersections and the problem raised by Curtis on an explicit determination of structure constants of $\calH$ 
are the motivation lying at the core of our work. The work is divided in two parts. Firstly, we determine an algorithm to decompose 
explicitly the intersection of Bruhat cells in terms of left $U$-coset representatives when $G$ is of small rank. Secondly, we use the output of this algorithm 
to determine the structure constants of $\calH$ in the case of adjoint types $\mathrm{A}_2$ and $\mathrm{B}_2$ and some generation %and relations 
properties in these algebras. 
The knowledge of the structure constants is essential for the determination of the $f_{T, \theta}$, see \cite[Section 4]{chang}. 
%HERE? The $f_{T, \theta}$ would certainly be obtained as in \cite{chang} with the use , as in
%of such structure constants. 
As the computations are involved and the methods need to be slightly adapted when moving from the simply connected 
case to the adjoint case, we postpone the problem of determining the $f_{T, \theta}$'s to subsequent work. %(HERE?)

We discuss now the methods and results of each of the two parts of the work. Let us focus on the algorithm described in the first part. %We recall that our 
The intersections of our interest have the form $UxU \cap zUy^{-1}U$ with $x, y$ and $z$ of the form $\dot{w}t_w$, with $\dot w$ a certain lift in $N_G(T)$ of an element $w$ of the Weyl group $W$ of $G$, and $t_w \in T$. %with $w$ some element of the Weyl group $W$ of $G$, $\dot w$ a certain lift of $w$ in $N_G(T)$, and $t_w \in T$. 
Further, $x$ comes with a fixed reduced expression in $W$. From results of Curtis \cite{curtis88,curtis09,curtis} which rely on work of Deodhar \cite{deodhar85}, we can decompose the intersection into disjoint sets %which are 
indexed by \emph{distinguished subexpressions} of the fixed reduced expression of $x$. Building on results of \cite{closure}, and exploiting the theory of the invariants for distinguished subexpressions, we first determine the intersections for certain $t_w$ and then we extend this result to all elements of $N_G(T)$. This decomposition method resembles the principle of ``collection from the left", see for example \cite[\S9.4]{Sim14}. This algorithm has been implemented in Python and is available in the GitHub repository \cite{github}. It has been successfully run for groups of small rank, in particular in types $\mathrm{A}_2$, $\mathrm{B}_2$ and $\mathrm{G}_2$ exhaustively. %Here is an overview of the main steps.
An outline of the algorithm is given in Section \ref{algo}.

We now move on to the investigation of the endomorphism algebras of Gelfand-Graev representations and their structure constants% of their endomorphism algebras
. We first mention some well-known results %in this direction
. The case $\bar G=\mathrm{SL}_2$ is worked out in detail in \cite{curtis93}. 
In \cite{chang}, the structure constants for $\calH$ in the case $\bar G=\GL_3$ are calculated and extensively used for the determination of the homomorphisms $f_{T, \theta}$. 
However, the knowledge of such structure constants does not give directly the ones for $\PGL_3(q)$. In fact, the endomorphism algebra with respect to $\PGL_3(q)$ is contained in the one with respect to $\GL_3(q)$, but \emph{not} a subquotient of it. In other words, 
a standard basis in the case of $\PGL_3(q)$ is not obtained from the one in 
the case of $\GL_3(q)$ in a natural way, and vice versa.%of $\calH'$ is not obtained in a natural way from the one of $\calH$. 
%Hence there is not a natural way to relate 

%This results in a difference of the parametrization of their standard bases. %used to determine the structure contants. %of the two algebras have 

%In this part of the work we %pursue this direction by investigating 
%We investigate 
%We focus in this part of the work on the investigation 
The focus of the second part of this work is on the description of the structure constants of endomorphism algebras of Gelfand-Graev representations for
%{\color{orange}
non-twisted groups %}
%non-twisted simple groups
of rank $2$ with connected center, by means of 
the determination of the $U$-coset representatives obtained in the first part 
which parametrize a standard basis for $\calH$ %a choice of a standard basis 
(described in \cite[Proposition 11.30]{CR81}, see also \S\ref{std_basis}). 
The choice of considering the \emph{adjoint} version of the simple algebraic group having $G$ as fixed-point subgroup %containing $G$ as fixed 
%
%fixed-point subgroups of simple adjoint algebraic groups 
(see \S\ref{setup_group}) is motivated by the following two facts. %: 1. 
Firstly, there is a unique %conjugacy class of characters when inducing characters in general position of a unipotent Sylow $p$-subgroup, 
Gelfand-Graev character, see for instance \cite[\S8.1]{carter_finite}. Secondly, it is easy to parametrize a maximally split torus which %makes the description of the standard basis of $\calH$ less cumbersome. 
allows a more compact description of the standard basis. %of $\calH$.
 
Without loss of generality, we may assume that a root datum for $G$ is chosen such that the character $\psi$ corresponding to a 
Gelfand-Graev representation restricts to the same character $\phi$ on each of the simple root groups. %parametrized by simple roots. 
%Our main result is as follows. 
We state our main result. 

\begin{theorem*} Let $q$ be a power of a prime $p$, and let $G=\PGL_3(q)$ for any $p$ or $G=\SO_5(q)$ for $p \ge 3$. Let $\psi$ be the Gelfand-Graev character of $G$. The structure constants of the endomorphism algebra of $\psi$  with respect to its standard basis are given in Table \ref{tab:A2} and Table \ref{tab:B2} respectively. In particular, the elements $e_1(a)$, $e_2(b)$ and $e_3$ for $a, b \in \FF_q^\times$ generate the endomorphism algebra of $\psi$ when $G=\PGL_3(q)$.
\end{theorem*}

The methods employed to obtain this result are as follows. %Our starting point is 
The formula given 
in \cite[Section 4]{curtis} (see \S\ref{std_basis}) 
%which 
gives a way to compute the structure constants of $\calH$ 
in terms of shifted $U$-double cosets intersections, which have been determined in the first part of our work for groups of small rank. %in rank $2$. 
%(REPET!) provided that one can express elements of $G$ in their canonical forms in terms of 
%a refinement of the Bruhat decomposition. (REPET!) 
The formula can then be written as a sum over certain %is a sum over certain 
indeterminates in $\FF_q$ and $\FF_q^\times$ of the character $\phi$ evaluated on 
a rational polynomial function in such indeterminates. %(REPET!) We first develop a general algorithm, which builds on 
%\cite{closure}, that determines the needed decomposition for the elements in $G$. This algorithm gives the desired expressions in low rank. (REPET!)
%, and can effectively be used when the ambient algebraic group has dimension less than 16. 
The remaining task is then to 
solve some equations over $\FF_q$ and to use these solutions, whose number %of solutions 
may vary with $q$, %and use these 
to express each sum coming from the formula in \cite[Section 4]{curtis} in a form which is as closed as possible. We can 
complete this step %in the cases 
for groups of type $A_2$ and $B_2$. 

In some cases, the structure constants 
involve certain generalizations of Gauss and Kloosterman sums (see \cite[Section 11]{IK}), and it seems not to be possible to 
further simplify those calculations. 
We provide the calculations explicitly in the most involved cases. Details for the remaining computations can be found on \cite{PS}. %difficult cases. The other
Although we were able to parametrize the standard basis of $\calH$ and determine the sums giving the structure constants when $G=\mathrm{G}_2(q)$, 
finding a satisfactory way of expressing these constants is still open. Here the major hurdle is the investigation of 
equations of degree $5$ or higher, 
whose generic sets of solutions are not easy to control.

We finish by presenting further directions for future work. An important open problem mentioned before is to 
describe explicitly the algebra homomorphisms $f_{T, \theta}: \calH \to \CC$. The methods are likely to 
involve the determination of the structure constants of $\calH$ and a minimal set of algebra generators, see 
\cite[Section 4--6]{chang}. The standard basis elements corresponding to the Weyl group elements $1$, $s_1s_2$ and $s_2s_1$ 
%are shown to 
generate $\calH$ in the cases of $\GL_3(q)$ and $\PGL_3(q)$. However, as explained in Section $\ref{constants}$, the techniques for the computations 
in type $\mathrm{A}_2$ do not generalize 
%even 
to type $\mathrm{B}_2$. We formulate the following question, an answer to which would generalize \cite[Theorem 2.1]{chang} to other types.
\begin{question*}
Can the algebra $\calH$ be always generated by elements of a standard basis which correspond to Weyl group elements of co-length less than or equal to $1$?
\end{question*}
\noindent Lastly, we mention the problem of finding a compact %generic 
way to express the structure constants of $\calH$, when the prime $p$ is large enough, for untwisted groups of type $\mathrm{G}_2$ or 
rank higher than $2$ and for twisted groups. %, and of generalizing 

The structure of the work is as follows. In 
%\S
Section \ref{pre} we describe the setup for %in more detail: 
%the 
groups, Gelfand-Graev characters and the standard basis of $\calH$, where structure constants of $\calH$ are expressed as sums over intersections of shifted $U$-double cosets. In Section \ref{algo} we show how these intersections can be obtained explicitly 
by means of the above mentioned algorithm. In %\S
Section \ref{sums} we discuss the sums which intervene in our description of the structure constants 
and we provide explicit computations %for the structure constants 
in the most involved case in type $B_2$. 
We collect in Section $\ref{constants}$ the structure constants in types $A_2$ and $B_2$.
\\

\noindent
\textbf{Acknowledgement:} The authors deeply thank G. Malle for his precious comments  
and feedback on an earlier version of the paper. Part of the work was developed 
during a research visit of the first author hosted at, and supported by, the Babe\c s-Bolyai University, 
and of the second author hosted at the Technische Universit\"at Kaiserslautern and 
supported by the SFB--TRR 195. The authors would like to thank both institutions 
for the kind hospitality.

\section{Preliminaries}
\label{pre}

\subsection{The group $G$}
\label{setup_group}
We consider a simple algebraic group $\bar G$ of adjoint type and $G=\bar G^F$, the fixed point subgroup under a Frobenius endomorphism $F$ which is not twisted. If $K$ is the subgroup of $G$ generated by unipotent elements, then with few exceptions, $K$ is simple \cite[Thm 2.2.7]{GLS}. By \cite[2.5.8(a)]{GLS} $G$ is an extension of $K$ by diagonal automorphisms. Moreover, since $\bar G$ is adjoint, we have $G=K\bar T^F$ for a split maximal torus $\bar T$ of $\bar G$ and $T=\bar T^{F}$ induces the full group of diagonal automorphisms of $K$.

The group $G$, viewed as extension of the finite simple group $K$ by all diagonal automorphisms, can be described in the Chevalley group setting as follows. We regard $K$ as an adjoint Chevalley group \cite[Theorem 11.1.2]{carter_simple} and consider the diagonal automorphisms described there. %with the description in \emph{loc.cit.}: 
The group $\bar T^{F}$ corresponds to $\hat H$ in \cite[\S7.1]{carter_simple}. %, see also \cite[p.200]{carter_simple} and \cite[Remark 2.5.11]{GLS} for characterizations of diagonal automorphisms.

Throughout, $T=\bar T^F\subseteq B=\bar B^F$ will denote a pair of a split maximal torus and Borel subgroup corresponding to fixed choices of the analogue $F$-stable subgroups $\bar T\subseteq \bar B$ in $\bar G$. The unipotent radical of $\bar B$ is $\bar U$ so $U$ is the Sylow $p$-subgroup of $B$. We fix the root system $\Phi$ with respect to $\bar T$. The simple roots $\Delta$ and the positive roots $\Phi^{+}$ are chosen with respect to $\bar B$. The root subgroups of $G$ are of the form $U_{\alpha}=\bar U_\alpha^{F}$ with $\bar U_{\alpha} \subseteq \bar B$ a root subgroup of $\bar G$. The subgroup of $U$ generated by root groups $U_\alpha$ corresponding to positive non-simple roots is called $U^{\ast}$.

Further, for $a, b \in \FF_q$ and $d \in \FF_q^\times$ we denote by $u_{\alpha}(a)$, $n_{\alpha}(b)$ the elements in $G$ which in \cite{carter_simple} are denoted by $x_{\alpha}(a)$, $n_{\alpha}(b)$ respectively and by $t_{\alpha}(d)$ the cocharacters corresponding to the coroots (see \cite[p.76]{carter_finite} or \cite[Remark 2.5.11(c)]{GLS}). In addition, $n_{\alpha}$ stands for $n_{\alpha}(1)$. To simplify notation, we write $i$ instead of $\alpha_i$ whenever $\alpha_i$ is an index, e.g. $u_{\alpha_1}(a)=u_{1}(a)$ and $t_{\alpha_6}(c)=t_{6}(c)$.

The root groups $U_{\alpha}$ are the images of $u_{\alpha}$ and the split maximal torus $T$ is generated by the elements of the form $t_{i}(c)$ %(In \cite{carter_simple}, this group is denoted by $\hat H$)
. %The $p$-Sylow subgroup generated by $U_\alpha$ for $\alpha\in\Phi^{+}$ is denoted by $U$.
A set of $\bar T$-coset representatives in the Weyl group $W=N_G(\bar T)/\bar T$ can be selected from words in the elements $n_{\alpha}\in G$ with $\alpha \in \Delta$; these choices will be made in the sequel. The set of these fixed representatives is $\dot W$ and for $w\in W$ we denote by $\dot w$ the corresponding representative. The simple reflections of $W$ are $\{s_{\alpha}=n_{\alpha}T:\alpha\in\Delta\}$, the length of a reduced word of $w\in W$ in terms of these simple reflections is denoted by $\ell(w)$ and the lifts of $s_{\alpha}$ are $\dot s_{\alpha}=n_{\alpha}\in\dot W$. In the sequel, we sometimes use the letter $\ell$ to denote a positive integer other than the length function; the meaning of $\ell$ is in any case determined by the local use. 

Using $t_{\alpha}$ to parametrize $T$ we have for every $\alpha, \beta \in \Phi^+$, $\lambda \in \FF_q^\times$ and $\mu \in \FF_q$ (see \cite[p.76]{carter_finite}),
\begin{equation}
  \label{eq:Torus}
        {}^{t_{\alpha}(\lambda)}u_{\beta}(\mu)=t_{\alpha}(\lambda)u_{\beta}(\mu)t_{\alpha}(\lambda)^{-1}=u_{\beta}(\lambda^{\delta_{\alpha,\beta}}\mu).
\end{equation}

The ground field of $\bar G$ is of characteristic $p>0$ and the Frobenius endomorphism $F$ is such that $|U_{\alpha}^{F}|=q$ for some power $q$ of $p$, i.e. $K\subseteq\bar G^{F}$ is a Chevalley group over the field $\FF_q$.
\subsection{Gelfand-Graev characters}\label{sub:GG} Following \cite{curtis}, let $\psi$ be a linear representation of $U$ and let $e=|U|^{-1}\sum_{u\in U}\psi(u^{-1})u\in\mathbb{C}G$ be the idempotent such that the induced representation $\psi^{G}$ is afforded by the module $\mathbb{C} Ge$. The Hecke algebra corresponding to $\psi$ is (as in \cite[\S11.D]{CR81})
$$
\calH=e\mathbb{C} Ge=\End_{\CC G}\CC Ge.
$$
Further, $\psi^{G}$ is a \emph{Gelfand-Graev character} if $\psi(U_{\alpha})\neq 1$ for simple roots $\alpha$ and $\psi(U_{\alpha})= 1$ for all other roots $\alpha$ \cite[\S8.1]{carter_finite}. %In this case we denote $H$ by $\mathcal{H}(G, U, \psi)$.
By \cite[Theorem 8.1.3]{carter_finite}, the algebra $\calH$ is proved to be abelian due to an isometry of Gelfand-Graev characters which has order two. 

The character $\psi$ has $U^{\ast}$ in its kernel. Since $U/U^{\ast}\cong \prod_{\alpha\in\Delta}U_{\alpha}\cong k^{|\Delta|}$, for any group isomorphism $u_{\alpha}:k\rightarrow U_{\alpha}$ we obtain characters $\phi_{\alpha}=(\psi|_{U_{\alpha}})\circ u_{\alpha}:k\rightarrow\CC$. Denoting by $\mathrm{Tr}:\FF_q \to \FF_p$ the field trace map, it then follows from \cite[Proposition 8.1.2]{carter_finite} that we may assume
\begin{equation}
  \label{eq:char_simple_roots}
\phi_{\alpha}=\phi\quad\forall\alpha\in\Delta, 
\qquad 
\text{ where }
\qquad 
\phi(x)=e^{\frac{2\pi i \mathrm{Tr}(x)}{p}}
.
\end{equation}
When determining the character values of $\psi$ and $\phi$, we will frequently use delta-notation for the description of the obtained structure constants, 
$$
\delta_{P}= %\left\{
\begin{cases}
  1 & \text{if $P$ is true,}\\
  0 & \text{if $P$ is false,}
  \end{cases}
%\right
\qquad
\text{ and }
\qquad
\delta_{a,b}=\delta_{a=b}= %\left\{
\begin{cases}
  1 & \text{if }a=b,\\
  0 & \text{if }a\neq b,
  \end{cases}
$$
for example $\delta_{1\in\{2,3\}}=0$. 

\subsection{Standard basis and structure constants}\label{std_basis} The set of $U$-double coset representatives of $G$ is $\{\dot wt:\,w\in W,\,t\in T\}$. The \emph{standard basis} $\calB$ of the algebra $\calH$ (with respect to $T\subseteq B$) is described in \cite[Proposition 11.30]{CR81}, it is
\begin{equation}
  \label{eq:Basis}
\calB=
\left\{
e_n=q^{\ell(w)}ene,
\quad\text{where}\quad
n=\dot wt,
\,w\in W, \,t\in T,
\quad\text{such that}\quad
{}^{n}\psi=\psi\text{ on }U\cap {}^{n}U
\right\}.
\end{equation}
For $e_\ell, e_m, e_n \in \mathcal{B}$ the corresponding structure constant is %(see \emph{loc.cit.})
\begin{equation}
\label{eq:GGR}
 [e_\ell e_m : e_n]%=\sum_{\mathbf{v} \in J(\mathbf{w}, x, y)} \sum_{(u, u_1, v) \in \mathcal{U}_{\mathbf{v}}}\psi((uu_1)^{-1}v)
 =\sum_{u\ell u'=nvm^{-1}\in U\ell U\cap nUm^{-1}U}\psi((uu')^{-1}v).
\end{equation}

An algorithm for obtaining these constants was provided by C. Curtis (see \cite{curtis}). A variation of this algorithm, with different choices for $U$-coset representatives was described in \cite{closure}. 
We use this latter version (see Section \ref{algo}) to obtain the elements
%$$
\begin{equation}
  \label{eq:doubleCells}
  u\ell u_1=nvm^{-1}\in U\ell U\cap nUm^{-1}U
  \end{equation}
%$$
explicitly in both their forms ($u\ell u_1$ and $nvm^{-1}$) for $G$ of type $A_2$ and $B_2$. 

Some of the more involved structure constants obtained here are expressed in terms of Gauss sums, Kloosterman sums and generalizations thereof %(these are 
described in Section \ref{sums}%)
. Similar calculations are obtained for $G_2$; however, finding concise formulae for the structure constants which  are meaningful in terms of character values of $\phi$ and symmetries of $\calH$ remains open in this case.

\subsection{Types $A_2$ and $B_2$}
We assume throughout that $\Phi$ is of type $A_2$ or $B_2$, with positive roots $\Phi^{+}$ ordered as follows,
$$
A_2:\,\alpha_1,\,\alpha_2,\,\alpha_1+\alpha_2, %\quad \text{for }A_2,
\qquad\text{and}\qquad
B_2:\,\alpha_1,\,\alpha_2,\,\alpha_1+\alpha_2,\,2\alpha_1+\alpha_2. %\quad \text{for }B_2
$$
The negative roots in $\Phi$ are ordered by $-\alpha_i=\alpha_{|\Phi^{+}|+i}$, e.g. $\alpha_6=-\alpha_1-\alpha_2$ in type $A_2$ and $\alpha_6=-\alpha_2$ in type $B_2$. The structure constants for the group $G$ are as in \cite[\S5.2]{carter_simple}. According to the choice of %with the 
extra-special pairs, %choosen such that in type $A_2$ we have
the non-trivial commutator relations among root elements are 
$$
{}^{u_{1}(x_{1})}u_{2}(x_{2})
=
u_{1}(x_{1})u_{2}(x_{2})u_{1}(x_{1})^{-1}
=
u_{2}(x_{2})u_{3}(x_{1}x_{2})
$$
in type $A_2$, and 
$$
{}^{u_{1}(x_{1})}u_{2}(x_{2})
=
u_{2}(x_{2})u_{3}(x_{1}x_{2})u_{4}(x_{1}^2x_{2})
\quad\text{and}\quad
{}^{u_{1}(x_{1})}u_{3}(x_{3})
=
u_{3}(x_{3})u_{4}(2x_{1}x_{3})
$$
in type $B_2$.

\subsubsection{}
\label{H_basis}
For the description of $\calB$, let $e_n\in\calB$ with $n=\dot wt$, $w\in W$ and $t\in T$. By \cite[(1) in proof of Thm. 49]{steinberg_lectures}, if $e_n\neq 0$ then  $w=w_0w_{\pi}$ where $w_0$ is the longest element of $W$ and $w_{\pi}$ is the longest element of the parabolic subgroup $W_{\pi}\leq W$ corresponding to the subset $\pi$ of simple roots. Notice that $\dot w$ and 
$\dot w_0\dot w_{\pi}$ might differ. 

In rank $2$, there are 4 standard parabolic subgroups of $W$. Their longest elements are $1$, $s_{1}$, $s_{2}$, $s_{1}s_{2}s_{1}$ for $A_2$ and $s_{1}s_{2}s_{1}s_{2}$ for $B_2$, so
$$
w_0w_{\pi}=\left\{
\begin{array}{ccccl}
 s_{1}s_{2}s_{1}, & s_{1}s_{2}, & s_{2}s_{1}, & 1 & \text{for } A_2,\\
 s_{1}s_{2}s_{1}s_{2}, & s_{2}s_{1}s_{2}, & s_{1}s_{2}s_{1}, & 1 & \text{for } B_2.\\
\end{array}
\right.
$$
%For the corresponding elements in $G$, we choose
We number the lifts in $N_G(T)$ of the above elements in the following way, 
$$
\begin{array}{cccccl}
\dot w_0=n_{1}n_{2}n_{1},& \dot w_1=n_{1}n_{2},& \dot w_2=n_{2}n_{1}, & \dot w_3=1 & \text{for } A_2,\\
\dot w_0=n_{1}n_{2}n_{1}n_{2},& \dot w_1=n_{2}n_{1}n_{2},& \dot w_2=n_{1}n_{2}n_{1}, & \dot w_3=1 & \text{for } B_2.\\
\end{array}
$$

\subsubsection{}\label{basis_rank_2}
Recall that the character ${}^{n}\psi$ is defined 
such that ${}^{n}\psi(u)=\psi(u^{n})$ for all 
$u \in U$. Thus in order to obtain the standard basis \eqref{eq:Basis}, %notice that ${}^{n}\psi(u)=\psi(u^{n})$ \cite[p.285]{CR81}, so, 
we need the condition %for elements in $\calB$ is
\begin{equation}
\label{basis_condition}
\psi(u)=\psi(u^{n})\quad\text{for}\quad u\in U\cap {}^{n}U
\end{equation}
to be satisfied. Hence
$$
\calB=
\left\{
e_n=q^{\ell(w)}ene,
\quad\text{where}\quad
n=\dot wt,
\quad
w=w_0w_{\pi}\in W,
\quad
t\in T,
\quad
\psi(u)=\psi(u^{n})
\text{ for }u\in U\cap {}^{n}U
\right\}.
$$
Using commutator relations, in particular \eqref{eq:Torus}, we can impose the condition $\psi(u)=\psi(u^{n})\text{ for }u\in U\cap {}^{n}U$ from \eqref{basis_condition} in our particular cases. Since $U\cap U^{\dot w_0}=1$, the condition for $w_0$ is empty, so $\dot w_0t\in\calB$ for all $t\in T$. For the element $n=n_1(a,b)=\dot w_1t_1(a)t_2(b)$, we have %$U\cap {}^{n}U=U_1$ and
$$
U\cap {}^{n}U=\left\{
\begin{array}{ll}
U_2
& \text{for } A_2,\\
U_1
& \text{for } B_2\\%\text{ and }G_2\\
\end{array}
\right.
\quad\text{and}\quad
\left\{
\begin{array}{ll}
  u_2(x)^{n_1(a,b)}=u_{1}(a^{-1}x)
& \text{for } A_2,\\
u_1(x)^{n_1(a,b)}=u_{1}(a^{-1}x)
& \text{for } B_2.\\
\end{array}
\right.
$$
With \eqref{eq:char_simple_roots}, we have that an element of the form $n_1(a,b)$ is in $\calB$ if and only if for all $x\in \FF_q$ we have %and $a,b\in \FF_q^{\times}$
\begin{equation}
  \label{w1_toral_condition}
  \left.
\begin{array}{r}
  \psi(u_2(x))=\psi(u_{1}(a^{-1}x))\\
    \psi(u_1(x))=\psi(u_{1}(a^{-1}x))
\end{array}
\right\}
  \Leftrightarrow
  \phi(x)=\phi(a^{-1}x)
  \Leftrightarrow a=1
 \text{ for } A_2\text{ and }B_2.\\
\end{equation}
Now consider the element $n_2(a,b)=\dot w_2t_1(a)t_2(b)$. We have %$U\cap {}^{n}U=U_2$ and
$$
U\cap {}^{n}U=\left\{
\begin{array}{ll}
U_1
& \text{for } A_2,\\
U_2
& \text{for } B_2\\%\text{ and }G_2\\
\end{array}
\right.
\quad\text{and}\quad
\left\{
\begin{array}{ll}
u_1(x)^{n_2(a,b)}=u_{2}(b^{-1}x)
& \text{for } A_2,\\
u_2(x)^{n_2(a,b)}=u_{2}(b^{-1}x)
& \text{for } B_2.\\
\end{array}
\right.
$$
Hence, as previously argued, we have $n_2(a,b)\in\calB$ if and only if $b=1$.

For $\dot w_3=1$ we consider $n=n_3(a,b)=t_1(a)t_2(b)$, so $U\cap {}^{n}U=U$ and $n_3(a,b)\in\calB$ if and only if
\begin{equation}
\label{w3_toral_condition}
\psi(u_1(x_1)u_2(x_2))=
%\left\{
%\begin{array}{ll}
  \psi(u_1(x_1a^{-1})u_2(x_2b^{-1}))
  \Leftrightarrow \phi(x_1+x_2)=\phi(a^{-1}x_1+b^{-1}x_2)
  \Leftrightarrow a=b=1
\end{equation}
in all cases. %Since we consider $p$ to be very good for $G$, condition \eqref{w3_toral_condition} implies that $a=b=1$.
These conditions determine $\calB$, namely, the elements in $\calB$ are parametrized by the elements
\begin{equation}
  \label{eq:nB}
n_0(a,b):=\dot w_0t_{1}(a)t_{2}(b),
\quad
n_1(c):=\dot w_1t_{2}(c),
\quad
n_2(d):=\dot w_2t_{1}(d),
\quad
n_3:=\dot w_3=1
%\quad\text{with }
%a,b,c,d\in\FF_q^{\times}.
\end{equation}
with $a,b,c,d\in\FF_q^{\times}$. 
The corresponding basis elements of $\calH$ are
$$
\calB=\big\{
e_0(a,b)=en_0(a,b)e,\quad
e_1(c)=en_1(c)e,\quad
e_2(d)=en_1(d)e,\quad
e_3=en_3e=e, 
\quad
a,b,c,d\in\FF_q^{\times}\big\}.
$$
We sometimes omit the parameters $a,b,c,d$ and write $e_i$ and $n_i$ for these elements.

\section{$U$-double cosets}
\label{algo}

\subsection{}\label{decomposition} Once the basis $\calB$ of $\calH$ is determined, we obtain the structure constants by \eqref{eq:GGR}. In order to do so, we need to determine all left $U$-coset representatives in the intersections \eqref{eq:doubleCells}:
%\begin{equation}
%\label{structure}
$$
Ux U\cap zUy^{-1}U\quad\text{for }x,y,z\in\calB.
$$
%\end{equation}
From results of C. Curtis \cite{curtis88,curtis09,curtis} building on results of V. Deodhar \cite{deodhar85}, the above intersection decomposes into disjoint sets parametrized by distinguished subexpressions of the fixed reduced expression for $x$,
$$
Ux U\cap zUy^{-1}U=\bigsqcup_{j\in J(x, y, z)}D_jU
\quad\text{for}\quad x,y,z\in N_G(T). 
$$
The sets $D_j$ of left $U$-coset representatives are described inductively. For concretely determining these intersections in both the $UxU$-form and $zUy^{-1}U$-form which are required for our purposes, one needs to translate the inductive descriptions into the two forms. By making use of the $D_j$ described in \cite{closure}, we illustrate how the algorithm works in general. %We recall that the algorithm is 
Such a procedure has been implemented in Python \cite{github}. 

The problem is divided in two steps. Since $x$, $y$ and $z$ are of the form $\dot wt_w$ for $w\in W$ and $t_w\in T$, the intersections can be obtained by
\begin{itemize}
\item[\textbf{I.}] determining first the intersections for certain $t_w$, and %when $z\in\dot W$ and for some $y,x$
\item[\textbf{II.}] deducing then the intersections for arbitrary $t_w$.
\end{itemize}
We show in \S\ref{example_B2} how the algorithm works in the case of $Un_0U\cap n_0Un_0U$ for $G$ of type $B_2$.

\subsection{Part I} \label{sub:pt1}The set $\dot W$ of reduced expressions in terms of simple reflections was chosen in \S\ref{H_basis}. Fix $x,y,z\in W$ and consider the fixed reduced expression $s_{i_n}\dots s_{i_1}$ of $x$, %which was fixed and 
which we may identify with the tuple $i=[i_n,\dots,i_1]$.
A subexpression $j\in\{0,i_n\}\times\dots\times \{0,i_1\}$ of $i$ has the following sets attached to it, 
$$
\begin{array}{rcl}
A_j&=&\{m:\ell(s_{j_m}\tau_{m-1}(j)y)<\ell(\tau_{m-1}(j)y)\},\\
B_j&=&\{m:\ell(s_{j_m}\tau_{m-1}(j)y)=\ell(\tau_{m-1}(j)y)\}=\{m:j_m=0\},\\
C_j&=&\{m:\ell(s_{j_m}\tau_{m-1}(j)y)>\ell(\tau_{m-1}(j)y)\},
\end{array}
$$
where $s_{0}=1$ and $\tau_{m}$ is the truncation map $j\mapsto\tau_{m}(j)=s_{j_m}\dots s_{j_1}$. %(with the convention that $s_0=1$).
Using these sets we define, for a tuple $\mu\in \FF_q^{\ell(x)}$,
$$
D_j(\mu):=\prod_{m=1}^{n} \big[u_{m}(\mu_m)\dot s_{m}\big]^{\delta_{m\in A_j}}\big[u_{-m}(\mu_m)\big]^{\delta_{m\in B_j}}\big[\dot s_{m}\big]^{\delta_{m\in C_j}}
$$
where the product is from right to left, with $\mu_m\in\FF_q$ if $m\in A_j$, $\mu_m\in \FF_q^{\times}$ if $m\in B_j$, and $\mu=1$ if $m\in C_j$ and where $\dot s_{\alpha}=n_{\alpha}$ and $u_{-m}=u_{-\alpha_m}$. Here the exponents are the delta-functions defined in \S\ref{sub:GG}. In this case, since $m$ is in exactly one of the sets $A_j$, $B_j$ or $C_j$, we take the term corresponding to the set in which $m$ lies.

We also define
$$
D_j'(\mu):=\prod_{m=1}^{n} \big[u_{m}(\mu_m)\dot s_{m}\big]^{\delta_{m\in A_j}}\big[u_{m}(\mu_m^{-1})\dot s_{m}(-\mu_m^{-1})u_{m}(\mu_m^{-1})\big]^{\delta_{m\in B_j}}\big[\dot s_{m}\big]^{\delta_{m\in C_j}}
$$
where the product is again from right to left and where $\mu_m$ is as before. We say that $m$ \emph{is of type} A\emph{,} B \emph{or} C according to the set $A_j$, $B_j$ or $C_j$ it lies in. It is a small calculation in type $A_1$ (using for example the relations in \cite{carter_simple}) to show that the $m$'s of type B are equal in the two expressions so we have
$$
D_j(\mu)=D_j'(\mu).
$$
It is also easy to notice that $D_j'(\mu)$ parametrizes distinct left $U$-coset representatives in $UxU$ (exactly $q^{|A_j|}(q-1)^{|B_j|}$ such representatives for the allowed values of $\mu$). Moreover, using commutator relations in \cite{carter_simple}, one can calculate the $UxU$-form of the element $D_j'(\mu)$, that is, 
$$
D_j'(\mu)=U_j'(\mu)xt_{\mu}V_j'(\mu),
$$
where $t_{\mu}\in T$ appears due to the B-type positions in $j$. Now, if $j$ is a distinguished subexpression of $i$, using \cite[Lemma 3.8]{closure} one can calculate the $zUy^{-1}U$-form of the element $D_j(\mu)$, namely 
$$
D_j(\mu)=zU_j(\mu)y^{-1}t_0V_j(\mu)
$$
where $1\neq t_0\in T$ may appear in the expression due to non-reduced products of elements in $\dot W$. 

It is in this part of the algorithm that the property of $j$ being a distinguished subexpression %for $j$ 
is used in an essential way. Namely this property allows one to move the root group elements corresponding to negative roots on the left-hand side of
$$
z^{-1}D_j(\mu)=U_j(\mu)y^{-1}t_0V_j(\mu)
$$
to the left, such that they end up in $U_j(\mu)\subseteq U$, (for details see \cite[\S3.3]{closure}). In this way we obtain
$$
D_j(\mu)U\subseteq Uxt_{\mu}U\cap zUy^{-1}t_0U.
$$
Since the elements $D_j(\mu)$ are distinct left $U$-coset representatives for different $j\in J(x,y,z)$ and %the 
different parameters $\mu$, it follows from \cite[Lemma 2.3]{curtis09} that we obtain all left $U$-cosets in $Uxt_{\mu}U\cap zUy^{-1}t_0U$.

\subsection{Part II}
We keep the assumptions and notation form Part I and consider the intersection
\begin{equation}
  \label{eq:algoII}
Uxt_xU\cap zt_zU(yt_y)^{-1}U
\end{equation}
for $t_x,t_y,t_z\in T$. Since
$
U_j'(\mu)xt_{\mu}V_j'(\mu)=D_j(\mu)=zU_j(\mu)y^{-1}t_0V_j(\mu)
$
we have
$$
U_j'(\mu)xt_{\mu}t_0\left[V_j'(\mu)V_j(\mu)^{-1}\right]^{t_0}=D_j(\mu)=zU_j(\mu)y^{-1}
$$
hence
$$
zt_zU_j(\mu)(yt_y)^{-1}=zt_zU_j(\mu)t_y^{-1}y^{-1}=
(t_z)^{z^{-1}}U_j'(\mu)xt_{\mu}t_0\left[V_j'(\mu)V_j(\mu)^{-1}\right]^{t_0}(t_y^{-1})^{y^{-1}}
$$
which equals
$$
\left[U_j'(\mu)\right]^{(t_z)^{z^{-1}}}\cdot x\cdot(t_z)^{z^{-1}x}t_{\mu}t_0(t_y^{-1})^{y^{-1}}\cdot\left[V_j'(\mu)V_j(\mu)^{-1}\right]^{t_0(t_y^{-1})^{y^{-1}}}.
$$
It follows that, for fixed $x,y,z,t_x,t_y,t_z$ the intersection
$$
Uxt_xU\cap zt_zU(yt_y)^{-1}U
$$
is non-empty if and only if
\begin{equation}
\label{toral_condition}
t_x=(t_z)^{z^{-1}x}t_{\mu}t_0(t_y^{-1})^{y^{-1}},
\end{equation}
in which case the elements
$$
\left\{
zt_zU_j(\mu)(yt_y)^{-1}
=\left[U_j'(\mu)\right]^{(t_z)^{z^{-1}}}xt_x\left[V_j'(\mu)V_j(\mu)^{-1}\right]^{t_0(t_y^{-1})^{y^{-1}}}
:\mu \in \FF_q^{\ell(x)},j \in J(x, y, z)\right\},
$$
where each entry $\mu_m$ is restricted to the conditions remarked at the beginning of \S\ref{sub:pt1}, form a complete set of left $U$-coset representatives for $Uxt_xU\cap zt_zU(yt_y)^{-1}U$. 

The elements in \eqref{eq:algoII} can thus be determined with the following steps:

\begin{figure}[ht]
  \centering
  \begin{minipage}{.95 \linewidth}
\begin{algorithm}[H]
\KwIn{Fixed reduced expressions for Weyl group elements. Elements $x,y,z\in W$ and $t_x,t_y,t_z\in T$.}
\KwResult{A set of left $U$-coset representatives in $U\dot xt_xU\cap \dot zt_z U (\dot yt_y)^{-1} U$.}
result $\leftarrow \emptyset$\;
\Begin{
   \ForEach{distinguished expression $j\in J(x,y,z)$}{
         consider the products $D_j(\mu)=D_j'(\mu)$\;
         bring $D_j(\mu)$ in the form $\dot zU_j(\mu)\dot y^{-1}t_0V_j(\mu)$\;
         bring $D_j'(\mu)$ in the form $U_j'(\mu)\dot xt_{\mu}V_j'(\mu)$\;
         \If{$t_x=(t_z)^{z^{-1}x}t_{\mu}t_0(t_y^{-1})^{y^{-1}}$}{
                 add $\dot zt_zU_j(\mu)(\dot yt_y)^{-1}=\left[U_j'(\mu)\right]^{(t_z)^{z^{-1}}}\dot xt_x\left[V_j'(\mu)V_j(\mu)^{-1}\right]^{t_0(t_y^{-1})^{y^{-1}}}$ to result\;
        }
 }
 \Return result\;
 }
% \caption{\"Uberblick}
% \label{algo}
\end{algorithm}
\end{minipage}
\end{figure}

In our calculations however, $t_x$, $t_y$ and $t_z$ do not range over $T$, namely they are restricted to the toral elements intervening in the standard basis %(see 
as in \eqref{eq:nB}%)
.

\subsection{$[e_0 e_0,e_0]$ in type $B_2$}
\label{example_B2}
%\iffalse
In order to determine the constants $[e_0(a_1,b_1)e_0(a_2,b_2):e_0(a_3,b_3)]$ in type $B_2$, we first need to determine the 
intersection in \eqref{eq:doubleCells}, namely 
$$
n_0(a_3,b_3)Un_0(a_2,b_2)^{-1}U\cap Un_0(a_1,b_1)U.
$$
The reduced expression for $w_0$ is $s_1s_2s_1s_2$ and the distinguished subexpressions in $J(w_0,w_0,w_0)$ are
$$
[0,0,0,0]\text{ of type BBBB, }
[0,2,0,2]\text{ of type BCBA and }
[1,0,1,0]\text{ of type CBAB.}
$$
\subsubsection{}\label{0202} For $j=[0,2,0,2]$ we have $D_j=u_{5}(x_1)n_{2}(1)u_{5}(x_3)u_{2}(x_4)n_{2}(1)$ so
$$
z^{-1}D_j=
\underbrace{u_{1}(-x_1)u_{2}(x_4)u_{3}(x_3)}_{U_j}
  \underbrace{n_{1}(-1)n_{2}(-1)n_{1}(-1)n_{2}(-1)}_{y^{-1}}\underbrace{t_{1}(-1)}_{t_0}.
$$
%where $t_0=t_1(-1)t_2(-1)$.
We also have $D_j'=u_{1}(x_1^{-1})n_{1}(-x_1^{-1})u_{1}(x_1^{-1})n_{2}(1)u_{1}(x_3^{-1})n_{1}(-x_3^{-1})u_{1}(x_3^{-1})u_{2}(x_4)n_{2}(1)=$
$$
%=
\underbrace{u_{1}(x_1^{-1})u_{2}\left(\frac{x_1(x_1x_4 - 2x_3)}{x_3^2}\right)u_{3}(x_3^{-1})}_{U_j'}
\underbrace{n_{1}(1)n_{2}(1)n_{1}(1)n_{2}(1)}_{=x}
\underbrace{t_{1}(x_1^2)t_{2}\left(\frac{x_3^2}{x_1^2}\right)}_{=t_\mu}
%\underbrace{t_{1}(x_1x_3)t_{2}(x_3^2)}_{=t_\mu}
\underbrace{u_{1}\left(\frac{x_3-x_1x_4}{x_1x_3}\right)u_{3}(x_3^{-1})u_{4}\left(\frac{x_4}{x_3^2}\right)}_{V_j'}
$$
which gives the first part of the algorithm. Condition \eqref{toral_condition} on the elements in $T$ is
$$
t_x=t_1(a_1)t_2(b_1)=t_{1}\left(-a_2a_3x_1^2\right)t_{2}\left(b_2b_3\frac{x_3^2}{x_1^2}\right),
\quad\text{where}\quad t_y=t_1(a_2)t_2(b_2)\quad\text{and}\quad t_z=t_1(a_3)t_2(b_3).
$$
Hence, for this distinguished expression the elements in $n_0(a_3,b_3)Un_0(a_2,b_2)^{-1}\cap Un_0(a_1,b_1)U$ are
$$
n_0(a_3,b_3)
u_{1}(-x_1)u_{2}(x_4)\cdots
n_0(a_2,b_2)^{-1}
=
u_{1}\left(\frac{1}{a_3x_1}\right)u_{2}\left(\frac{x_1(x_1x_4 - 2x_3)}{b_3x_3^2}\right)\cdots
n_0(a,b)
u_{1}\left(\frac{x_1x_4-x_3}{a_2x_1x_3}\right)\cdots
$$
where $\mu=(x_1,1,x_3,x_4)\in \FF_q^{\times}\times\{1\}\times \FF_q^{\times}\times \FF_q$ and where the dots indicate terms in $U_3$ and $U_4$. For our purposes, we are only interested in $U_1$ and $U_2$, so, once $D_j=D_j'$ has been brought in the %two forms, 
$zUy^{-1}U$-form and in the $UxU$-form, we may ignore the terms in $U_3$ and $U_4$, in that they lie in $\ker\psi$ and do not contribute to the sum \eqref{eq:GGR} with distinct summands. 

\subsubsection{}\label{1010} For $j=[1,0,1,0]$ we have $D_j=n_{1}(1)u_{6}(x_2)u_{1}(x_3)n_{1}(1)u_{6}(x_4)$ so
$$
z^{-1}D_j=
\underbrace{u_{1}(x_3)u_{2}(-x_4)u_{4}(-x_2)}_{U_j}
\underbrace{n_{1}(-1)n_{2}(-1)n_{1}(-1)n_{2}(-1)}_{y^{-1}}.
$$
Here $t_0=1$. We also have $D_j'=n_{1}(1)u_{2}(x_2^{-1})n_{2}(-x_2^{-1})u_{2}(x_2^{-1})u_{1}(x_3)n_{1}(1)u_{2}(x_4^{-1})n_{2}(-x_4^{-1})u_{2}(x_4^{-1})=$
$$
=
\underbrace{u_{2}\left(\frac{1}{x_4} + \frac{x_3^2}{x_2}\right)u_{3}\left(-\frac{x_3}{x_2}\right)u_{4}\left(x_2^{-1}\right)}_{U_j'}
\underbrace{n_{1}(1)n_{2}(1)n_{1}(1)n_{2}(1)}_{=x}
\underbrace{t_{1}\left(\frac{x_2}{x_4}\right)t_{2}\left(x_4^2\right)}_{=t_\mu}
%\underbrace{t_{1}(-x_2)t_{2}(x_2x_4)}_{=t_\mu}
\underbrace{u_{1}\left(\frac{x_3x_4}{x_2}\right)u_{2}\left(x_4^{-1}\right)u_{4}\left(x_2^{-1}\right)}_{V_j'}.
$$
The condition \eqref{toral_condition} on elements in $T$ is
$$
t_x=t_1(a_1)t_2(b_2)=t_{1}\left(a_2a_3\frac{x_2}{x_4}\right)t_{2}\left(b_2b_3x_4^2\right),
\quad\text{where}\quad t_y=t_1(a_2)t_2(b_2)\quad\text{and}\quad t_z=t_1(a_3)t_2(b_3).
$$
For this distinguished expression the elements in $n_0(a_3,b_3)Un_0(a_2,b_2)^{-1}U\cap Un_0(a_1,b_1)$ are
$$
n_0(a_3,b_3)
u_{1}(x_3)u_{2}(-x_4)\cdots
n_0(a_2,b_2)^{-1}
u_{2}\left(\frac{x_2 + x_3^2x_4}{b_3^2x_2x_4}\right)\cdots
=
u_{1}\left(\frac{x_3x_4}{a_2x_2}\right)u_{2}\left(\frac{1}{b_2x_4}\right)\cdots
n_0(a,b)
$$
where $\mu=(1,x_2,x_3,x_4)\in \{1\}\times \FF_q^{\times}\times \FF_q\times \FF_q^{\times}$ and where the dots indicate terms in $U_3$ and $U_4$.
%\fi

%
\subsubsection{}\label{0000}For $j=[0,0,0,0]$ we have $D_j=u_{5}(x_1)u_{6}(x_2)u_{5}(x_3)u_{6}(x_4)=u_{-1}(x_1)u_{-2}(x_2)u_{-1}(x_3)u_{-2}(x_4)$, so
$$
z^{-1}D_j=
  \underbrace{u_{1}(-x_1 - x_3)u_{2}(-x_2 - x_4)u_{3}(-x_2x_3)u_{4}(-x_2x_3^2)
    }_{U_j}
    \underbrace{n_{1}(-1)n_{2}(-1)n_{1}(-1)n_{2}(-1)}_{=(n_1n_2n_1n_2)^{-1}=y^{-1}}.
$$
%where $t_0=t_1(-1)$.
Moreover, %Further,
$$
D'_j=u_{1}(x_1^{-1})n_{1}(-x_1^{-1})u_{1}(x_1^{-1})u_{2}(x_2^{-1})n_{2}(-x_2^{-1})u_{2}(x_2^{-1})u_{1}(x_3^{-1})n_{1}(-x_3^{-1})u_{1}(x_3^{-1})u_{2}(x_4^{-1})n_{2}(-x_4^{-1})u_{2}(x_4^{-1})
$$
and calculating the Bruhat form, one finds that
$$
D'_j=
\underbrace{u_{1}\left(x_1^{-1}\right)u_{2}\left(\frac{x_1^2x_2 + x_4(x_1(x_1 + 2x_3) + x_3^2)}{x_2x_3^2x_4}\right)u_{3}\left(-\frac{x_1 + x_3}{x_1x_2x_3}\right)u_{4}\left(\frac{1}{x_1^2x_2}\right)
  }_{=U_j'}
\underbrace{n_{1}(1)n_{2}(1)n_{1}(1)n_{2}(1)}_{=x}\cdot
$$
$$
\cdot \underbrace{t_{1}\left(\frac{x_{1}^2x_{2}}{x_{4}}\right)t_{2}\left(\frac{x_{3}^2x_{4}^2}{x_{1}^2}\right)}_{t_\mu}
%\cdot \underbrace{t_{1}(-x_1x_2x_3)t_{2}(x_2x_3^2x_4)}_{t_\mu}
\underbrace{u_{1}\left(\frac{x_1x_2 + x_4(x_1 + x_3)}{x_1x_2x_3}\right)u_{2}\left(\frac{1}{x_4}\right)u_{3}\left(-\frac{1}{x_3x_4}\right)u_{4}\left(\frac{x_2 + x_4}{x_2x_3^2x_4}\right)}_{=V_j'}.
$$
For the second part of the algorithm we have $t_x=t_1(a_1)t_2(b_1)$, $t_y=t_1(a_2)t_2(b_2)$, $t_z=t_1(a_3)t_2(b_3)$ and condition \eqref{toral_condition} on elements in $T$ is
$$
t_1(a_1)t_2(b_1)
=t_x
=(t_z)^{z^{-1}x}t_{\mu}t_0(t_y^{-1})^{y^{-1}}
=t_{1}\left(\frac{a_{2}a_{3}x_{1}^2x_{2}}{x_{4}}\right)t_{2}\left(\frac{b_{2}b_{3}x_{3}^2x_{4}^2}{x_{1}^2}\right). 
%=t_{1}(-a'a''x_1x_2x_3)t_{2}(b'b''x_2x_3^2x_4)
$$
In fact, here $z^{-1}x=1$ and $t_y^{y^{-1}}=t_y^{-1}$. Continuing, the elements in $n_0(a_3,b_3)Un_0(a_2,b_2)^{-1}U\cap Un_0(a,b)$ are
%corresponding to this distinguished expression $j$ are
$$
n_0(a_3,b_3)
\underbrace{u_{1}(-x_1 - x_3)u_{2}(-x_2 - x_4)\cdots}_{U_j(\mu)}
%u_{1}(x_{1})u_{2}(x_{2})..
n_0(a_2,b_2)^{-1}
=
$$
$$
\underbrace{
  u_{1}\left(\frac{1}{a_{3}x_{1}}\right)u_{2}\left(\frac{(x_{1}^2x_{2} + x_{4}(x_{1}(x_{1} + 2x_{3}) + x_{3}^2))}{b_{3}x_{2}x_{3}^2x_{4}}\right)
%  u_{1}\left(\frac{b''}{(a'')^2x_1}\right)u_{2}\left(\frac{(a'')^2(x_1^2x_2 + x_1^2x_4 + 2x_1x_3x_4 + x_3^2x_4}{(b'')^2{x_2}x_3^2x_4}\right)..
}_{
  %\left[V_j'(\mu)^{-1}V_j(\mu)\right]^{t_0(t_y^{-1})^{y^{-1}}}=\left[V_j'(\mu)^{-1}\right]^{(t_y^{-1})^{y^{-1}}}
  \left[U_j'(\mu)\right]^{(t_z)^{z^{-1}}}
}
\cdots
%\cdot\underbrace{u_{1}\left(\frac{b'(x_1x_2 + x_1x_4 + x_3x_4)}{(a')^2x_1x_2x_3}\right)u_{2}\left(\frac{(a')^2}{(b')^2x_4}\right)
%..}_{U_j'(\mu)^{(t_z)^{z^{-1}}}}
n_0(a_1,b_1)
\underbrace{
  u_{1}\left(\frac{x_{1}x_{2} + x_{4}(x_{1} + x_{3})}{a_{2}x_{1}x_{2}x_{3}}\right)u_{2}\left(\frac{1}{b_{2}x_{4}}\right)\cdots
}
  _{
    \left[V_j'(\mu)V_j(\mu)^{-1}\right]^{t_0(t_y^{-1})^{y^{-1}}}=V_j'(\mu)^{(t_y^{-1})^{y^{-1}}}
    }
$$
where $\mu=(x_1,x_2,x_3,x_4)\in (\FF_q^{\times})^4$, where the dots indicate terms in the root groups $U_3$ and $U_4$ and where $t_x$, $t_y$ and $t_z$ satisfy the relation previously deduced.

\section{Calculating structure constants}
\label{sums}
We describe in this section the methods employed to obtain 
the values of the structure constants of endomorphism algebras of 
Gelfand-Graev representations in types $A_2$ and $B_2$. Such methods involve manipulating 
some well-known sums over $\FF_q$ %and their properties. These 
which generalize the so-called 
\emph{Gauss sums} over finite fields and involve values of the character 
$\phi$ fixed in \eqref{eq:char_simple_roots}. %and some of their properties, . 
These sums, 
in particular quadratic Gauss sums and 
Kloosterman sums, 
frequently occur in the description of the 
structure constants in Tables \ref{tab:A2} and \ref{tab:B2} and allow a 
more compact way of collecting the structure constants.  

\subsection{} \label{sub:firG}We first remark some properties of roots in finite fields. 
Let $\mathbb{F}_q^\times=\langle \mu \rangle$, and let $r$ be a prime. 
If $r \nmid (q-1)$ then the map $x \mapsto x^r$ is an isomorphism, and every element of $\mathbb{F}_q^\times$ has a 
unique $r$-th root. 
Let us now assume $r \mid (q-1)$, and let $r'$ be such that $rr'=q-1$. In this case, let
$$\FF_{q, r}^\times:=\langle \mu^{r} \rangle=\{x \in \FF_q^\times \mid x = y^r \text{ for some }y \in \FF_q^\times\}.$$ 
Notice that the map $x \mapsto x^r$ has kernel of size $r$, hence its image $\FF_{q, r}^\times$ has cardinality 
$r'$. An element 
has an $r$-th root in $\mathbb{F}_q$ if and only if it is one of the $r'$ elements 
in the image of $x \mapsto x^r$. Moreover, we have $(\mu^{r'})^r=1 \in \mathrm{im}(x \mapsto x^r)$. 
Thus %if $B \in \mathbb{F}_q$ is an element having an $r$-th root, 
%then it 
each $B \in \FF_{q,r}^\times$ has exactly $r$ distinct $r$-th roots in $\mathbb{F}_q$. 

We recall the following properties obtained 
by substitutions in some Gauss sums. Let $a, b\in \FF_q$ and $c\in \FF_q^\times$%, and let $f:\FF_q \to \CC$
. Then we have %, see also \cite[Lemma 1.2]{chang}. 
%\begin{lemma}\label{lem:Chang}
\begin{itemize}
\item[(i)] $\sum_{x \in \FF_q}\phi(x)=0$ and $\sum_{x \in \FF_q^\times}\phi(ax)=-1+q\delta_{a, 0}$.
\item[(ii)] $\sum_{x \in \FF_q^\times}\phi(ax+\frac{b}{x})
=\sum_{x \in \FF_q^\times}\phi(x+\frac{ab}{x})+q\delta_{a, 0}\delta_{b, 0}$.
\item[(iii)] If $d \in \FF_q^\times$ has $r$ distinct roots in $\FF_q^\times$, then $\sum_{\zeta^r=d}\sum_{x \in \FF_q}
\phi(x(a-c\zeta))
%+f(\zeta))
=q\delta_{d, (\frac{a}{c})^r}
%\phi(f(\frac{a}{c}))
$. 
\end{itemize}
For (i) and (ii), we %recall that $\sum_{x \in \FF_q}\phi(x)=0$ and we 
refer to \cite[Lemma 1.2]{chang}. We now show (iii). Denote by 
$\zeta_1, \dots, \zeta_r$ the $r$ distinct $r$-th roots of $d$. 
If $d \ne (a/c)^r$ then $a-c\zeta_i \ne 0$ for every $i=1, \dots, r$, while 
if $d=(a/c)^r$ then $a-c\zeta_{\hat{i}} = 0$ for exactly one index $\hat{i} \in \{1, \dots, r\}$. 
The claim now follows by applying (i). 

\subsection{}\label{sub:Gauss}
The \emph{standard quadratic Gauss sum} has the form 
\begin{equation}\label{eq:Gauss}
G:=\sum_{x \in \mathbb{F}_q} \phi(x^2).
\end{equation}
We first assume that $p$ is an odd prime. %from 
Let $a \in \mathbb{F}_q^\times$. Then  it is easy to see that  
$$\sum_{x \in \mathbb{F}_q} \phi(ax^2)=
%\begin{cases}
(2\delta_{a \in \FF_{q, 2}^\times}-1)G, %, & \text{ if } a \text{ is a square in } \mathbb{F}_q, \\
%-G, & \text{ otherwise; }
%\end{cases}
$$
for example by applying (i) in \S\ref{sub:firG}. For $q=p$, we have that this sum evaluates to $\sqrt{p}$ if $p$ is congruent to $1$ 
modulo $4$, and to $i\sqrt{p}$ if $p$ is congruent to $-1$ modulo $4$. For generic odd $q$, 
we refer to \cite[Chapters 3 and 11]{IK}; 
%one just knows that $|G|=\sqrt{q}$, but an 
in particular, an explicit value of this sum is in general not known.  %(and possibly does not exist). 

It is then straightforward to obtain the following generalization,
$$
\sum_{x \in \mathbb{F}_q} \phi(Ax^2+Bx+C)=
%\begin{cases}
(2\delta_{a \in \FF_{q, 2}^\times}-1)G\phi\left(C-\frac{B^2}{4A}\right). 
$$

Quadratic gauss sums and their generalizations yield different sum values if $p=2$. %As remarked in the 
%Introduction, this is one of the reasons why we treat the case of $\SO_5(q)$ just for $p \ge 3$. 
Namely in this case %if $p=2$ then %fact Notice that if $p = 2$ then 
%this case 
we have that $G$ as defined in \eqref{eq:Gauss} evaluates to $0$, as 
the map $x \mapsto x^2$ is an automorphism of 
$\mathbb{F}_q$ when $q=2^f$. Moreover, we recall that in this case we have that 
$\phi(x^2+x)=1$ for every $x \in \mathbb{F}_q$. %and that 
%Writing %$A$ is always a square in $\mathbb{F}_q$, 
%$A=a^2$ for each $A \in \FF_q^\times$ $a \in \mathbb{F}_q$. 
As a consequence, if we let $A, B, C \in \FF_q^\times$, then we have %by writing $A=a^2$ for some $a \in \FF_q^\times$ we get 
\begin{align*}
\sum_{x \in \mathbb{F}_q} \phi(Ax^2+Bx+C)=
\begin{cases}
q\phi(C), 
& \text{ if } A=B^2 \\
0, 
& \text{ otherwise. }
\end{cases}
\end{align*}

\subsection{} Let $\ell \mid q-1$. We define the following sum, for $B, a, b, a', b' \in \FF_q$,
$$\tilde{\calS}_\ell(B, a, b, a', b'):=\sum_{\zeta \in \FF_q^\times \mid \zeta^\ell =B}
\phi\left(a'\zeta^2+a\zeta+\frac{b}{\zeta}+\frac{b'}{\zeta^2}\right).
$$
A special case of the sum $\tilde{\calS}_\ell(B, a, b, a', b')$ repeatedly appears in Tables \ref{tab:A2} and \ref{tab:B2}, namely
$$\calS_\ell(B, a, b):=\tilde{\calS}_\ell(B, a, b, 0, 0)=\sum_{\zeta \in \FF_q^\times \mid \zeta^\ell =B}
\phi\left(a\zeta+\frac{b}{\zeta}\right).
$$
We call $\tilde{\calS}_\ell(B, a, b, a', b')$ and $\calS_\ell(B, a, b)$ \emph{generalized Kloosterman 
sums}, as the value of $-\calS_{q-1}(1, a, b)$ is %given as 
equal to the value of the Kloosterman sum in \cite[\S11.5]{IK} when the maps $\psi$ and $\varphi$ appearing there 
%equal to 
are multiplications by $a$ and $b$ respectively. Notice that $\tilde{\calS}_\ell(1, a, b, a', b')
=\tilde{\calS}_\ell(1, b, a, b', a')$, thus $\calS_{\ell}(1, a, b)=\calS_{\ell}(1, b, a)$. Although 
$|\calS_{\ell}(1, a, b)|$ can be sharply bounded, no simpler formula for the sum $\calS_\ell(B, a, b)$ exists, 
see \cite[\S1.4]{IK}.

\subsection{$[e_0e_0:e_0]$ in type $B_2$}\label{sub:S000}
We are now in a position to continue the calculation of the structure constants with the use of \eqref{eq:GGR}. 
We provide here in full details the computations giving the structure constants 
$[e_0(a_1, b_1)e_0(a_2, b_2):e_0(a_3, b_3)]$ in type $B_2$. %All the other structure constants 
%obtained in 
This is the most complicated case to study from Tables \ref{tab:A2} and \ref{tab:B2}. The computational details for the other structure constants 
in such tables are collected in \cite{PS} and can be obtained more easily by using the same methods. %can be obtained more easily by the same methods. %in an easier 
%

%in what follows the most difficult example
%Because 
Since there are three distinguished subexpressions $j$, the sum describing these structure constants splits into three parts which we treat separately.

\subsubsection{$j=[0, 2, 0, 2]$} By \S\ref{0202}
we have that
\begin{equation}
\label{toral0202}
a_{1}=-a_{2}a_{3}x_{1}^2\quad\text{and}\quad
b_{1}=\frac{b_{2}b_{3}x_{3}^2}{x_{1}^2}.
\end{equation}
The subsum in \eqref{eq:GGR} restricted to the elements in \S\ref{0202} is
$$
\sum\psi\left(u_{1}\left(-x_{1} - \frac{1}{a_{3}x_{1}} - \frac{x_{4}}{a_{2}x_{3}} + \frac{1}{a_{2}x_{1}}\right)u_{2}\left(x_{4} - \frac{x_{1}^2x_{4}}{b_{3}x_{3}^2} + \frac{2x_{1}}{b_{3}x_{3}}\right)\right)$$
with summation over $(x_{1},x_{2},x_{3},x_4)\in \FF_q^{\times}\times\{1\}\times \FF_q^{\times}\times \FF_q$. This equals %$x_{2}=1$, $x_{3}\in k^{\ast}$, $x_{4}\in k$
$$
\sum\phi\left(-x_{1} - \frac{1}{a_{3}x_{1}} - \frac{x_{4}}{a_{2}x_{3}} + \frac{1}{a_{2}x_{1}}+x_{4} - \frac{x_{1}^2x_{4}}{b_{3}x_{3}^2} + \frac{2x_{1}}{b_{3}x_{3}}\right).
$$
%\textcolor{red}{Here assume $q$ odd!} Call $S$ the set of squares in $k^\times$. 
%Assume that 
If $-\frac{a_1}{a_2a_3} \notin \FF_{q, 2}^\times$ or $\frac{b_1}{b_2b_3} \notin \FF_{q, 2}^\times$, 
then the equations \eqref{toral0202} do not have a common solution, and the sum is zero. 
We assume $-\frac{a_1}{a_2a_3}, \frac{b_1}{b_2b_3}\in \FF_{q, 2}^\times$. %and $$ both 
%have a square root, say $\alpha$ and $\beta$ respectively. 
Then we have that 
$$(x_1, x_3) \in \left\{(\zeta_1, \zeta_1\zeta_2) \mid \zeta_1^2=-\frac{a_1}{a_2a_3}, \zeta_2^2=\frac{b_1}{b_2b_3}\right\}.$$
Hence we can write the 
sum as 
\begin{align*}
\sum_{ \zeta_1^2=-\frac{a_1}{a_2a_3}}
%\left(
\sum_{\zeta_2^2=\frac{b_1}{b_2b_3}}
\sum_{x_4 \in k} 
\phi(
x_4(1-\frac{b_2}{b_1}+\frac{a_3b_2b_3\zeta_1}{a_1b_1}\zeta_2)
-\zeta_1+\frac{a_2}{a_1}\zeta_1-\frac{a_3}{a_1}\zeta_1+2\frac{b_2}{b_1}\zeta_2
),
\end{align*}
and by (iii) %of Lemma \ref{lem:Chang} 
in \S\ref{sub:firG} this can be written as 
$$
\sum_{ \zeta_1^2=-\frac{a_1}{a_2a_3}}
\sum_{\zeta_2^2=\frac{b_1}{b_2b_3}}
q\delta_{\frac{b_2}{b_1}-1, \frac{a_3b_2b_3}{a_1b_1}\zeta_1\zeta_2}
\phi
\left(
-\zeta_1+\frac{a_2}{a_1}\zeta_1-\frac{a_3}{a_1}\zeta_1+2\frac{b_2}{b_1}\zeta_2
\right)
.
$$

\iffalse
\begin{align*}
\sum_{ \zeta_1^2=-\frac{a_1}{a_2a_3}}
q\delta_{\frac{b_1}{b_2b_3}, -\frac{a_1a_2(b_1-b_2)^2}{a_3b_2^2b_3^2}}
\phi\left(-\zeta_1+\frac{a_2}{a_1}\zeta_1-\frac{a_3}{a_1}\zeta_1-2\frac{b_2}{b_1}(1-\frac{b_2}{b_1})\frac{a_1b_1}{a_3b_2b_3\zeta_1}\right). 
\end{align*}
This sum evaluates
%This sum, which we call $S_{0,0}^{0, 2}$ in Table \ref{tab:B2}, evaluates
to zero if $b_1 =b_2$, and to 
%which is zero if $b_1 =b_2$, and %Otherwise, we can write it as 
\begin{align*}
q\delta_{b_3, -\frac{(b_1-b_2)^2a_1a_2}{a_3b_1b_2}}\sum_{ \zeta_1^2=-\frac{a_1}{a_2a_3}}
\phi\left(-\zeta_1+\frac{a_2}{a_1}\zeta_1-\frac{a_3}{a_1}\zeta_1+2\frac{a_3b_2}{a_1(b_1-b_2)}\zeta_1
%prev: -2\frac{a_2(b_1-b_2)}{b_1b_3}\zeta_1
\right)
\end{align*}
otherwise. 
\fi
\subsubsection{$j=[1, 0, 1, 0]$} By \S\ref{1010} we have that
$$
a_1=\frac{a_{2}a_{3}x_{2}}{x_{4}}
\quad\text{and}\quad
b_1=b_{2}b_{3}x_{4}^2
.
$$
The subsum in \eqref{eq:GGR} restricted to the elements in \S\ref{1010} is
\begin{align*}
\sum\psi & \left(
u_{1}\left(x_{3} - \frac{x_{3}x_{4}}{a_{2}x_{2}}\right)
u_{2}\left(-x_{4} - \frac{1}{b_{3}x_{4}} - \frac{x_{3}^2}{b_{3}x_{2}} - \frac{1}{b_{2}x_{4}}\right)\right)\\
=&\sum\phi\left(
x_{3} - \frac{x_{3}x_{4}}{a_{2}x_{2}}
-x_{4} - \frac{1}{b_{3}x_{4}} - \frac{x_{3}^2}{b_{3}x_{2}} - \frac{1}{b_{2}x_{4}}\right),
\end{align*}
with summation over $(x_{1},x_{2},x_{3},x_4)\in \{1\}\times \FF_q^{\times}\times \FF_q\times \FF_q^{\times}$. 

If $\frac{b_1}{b_2b_3} \notin \FF_{q, 2}^\times$, then we have no 
possible value for $x_4$, and the sum is zero. 
We then assume $\frac{b_1}{b_2b_3} \in \FF_{q, 2}^\times$. %, and denote by $\beta$ a square root of $b_1/(b_2b_3)$. 
In this case, we have 
$$(x_2, x_4) \in \left\{(\frac{a_1}{a_2a_3}\zeta, \zeta) \mid \zeta^2=\frac{b_1}{b_2b_3}\right\},$$ 
and we can write the above sum as 
\begin{align*}
%&
\sum_{\zeta^2=\frac{b_1}{b_2b_3}}\sum_{x_3 \in k} \phi(x_3-\frac{x_3}{a_2}\frac{\zeta a_2 a_3}{a_1 \zeta}-\zeta-\frac{1}{b_3\zeta}-\frac{x_3^2}{b_3}\frac{a_2a_3}{a_1\zeta}-\frac{1}{b_2\zeta})%+\sum_{x_3 \in k} \phi(x_3-\frac{x_3}{a_2}\frac{\zeta a_2 a_3}{a_1 \zeta}+\zeta+\frac{1}{b_3\zeta}+\frac{x_3^2}{b_3}\frac{a_2a_3}{a_1\zeta}+\frac{1}{b_2\zeta})
=%&
\sum_{\zeta^2=\frac{b_1}{b_2b_3}}\sum_{x_3 \in k}\phi\left(A(\zeta)x_3^2+Bx_3+C(\zeta)\right)%+\sum_{x_3 \in k}(-Ax_3^2+Bx_3-C)
,
\end{align*}
where $A(\zeta)=\frac{a_2a_3}{a_1b_3\zeta}$, $B=1-\frac{a_3}{a_1}$ and $C(\zeta)=\zeta+\frac{1}{b_2\zeta}+\frac{1}{b_3\zeta}$. 
The above sum can
%, which we call $S_{0, 0}^{0, 1}$ in Table \ref{tab:B2}, can
be written by \S\ref{sub:Gauss} as
%By \S\ref{sub:Gauss}, the sum can be written as 
$$G\sum_{\zeta^2=\frac{b_1}{b_2b_3}}(2\delta_{A(\zeta) \in \FF_{q,2}^\times}-1)\phi\left(C(\zeta)-\frac{B^2}{4A(\zeta)}\right).$$
%$$(-1)^{\delta_{A \notin S}}G\phi(C-\frac{B^2}{4A})+(-1)^{\delta_{-A \notin S}}G\phi(-C+\frac{B^2}{4A}).$$

%
\subsubsection{$j=[0, 0, 0, 0]$} 
In this case, the subsum in \eqref{eq:GGR} restricted to the elements in \S\ref{0000} is 
$$\sum\psi(u_{1}(-x_{1} - x_{3} - \frac{1}{a_{3}x_{1}} - \frac{1}{a_{2}x_{3}} - \frac{x_{4}}{a_{2}x_{2}x_{3}} - \frac{x_{4}}{a_{2}x_{1}x_{2}})u_{2}(-x_{2} - x_{4} - \frac{x_{1}^2}{b_{3}x_{3}^2x_{4}} - \frac{x_{1}^2}{b_{3}x_{2}x_{3}^2} - \frac{2x_{1}}{b_{3}x_{2}x_{3}} - \frac{1}{b_{3}x_{2}} - \frac{1}{b_{2}x_{4}}))$$
$$
=\sum\phi(
-x_{1} - x_{3} - \frac{1}{a_{3}x_{1}} - \frac{1}{a_{2}x_{3}} - \frac{x_{4}}{a_{2}x_{2}x_{3}} - \frac{x_{4}}{a_{2}x_{1}x_{2}}
-x_{2} - x_{4} - \frac{x_{1}^2}{b_{3}x_{3}^2x_{4}} - \frac{x_{1}^2}{b_{3}x_{2}x_{3}^2} - \frac{2x_{1}}{b_{3}x_{2}x_{3}} - \frac{1}{b_{3}x_{2}} - \frac{1}{b_{2}x_{4}}
),
$$
and by \S\ref{0000} we have %that
%(of type BBBB) $x_{1}\in k^{\ast}$, $x_{2}\in k^{\ast}$, $x_{3}\in k^{\ast}$, $x_{4}\in k^{\ast}$, :
%$$t_{1}(a_{1})t_{2}(b_{1})=t_{1}(a_{2}a_{3}x_{1}^2x_{2}/x_{4})t_{2}(b_{2}b_{3}x_{3}^2x_{4}^2/x_{1}^2)$$
%We have that
$a_{1}=a_{2}a_{3}x_{1}^2x_{2}/x_{4}$ 
and 
$b_{1}=b_{2}b_{3}x_{3}^2x_{4}^2/x_{1}^2$.
The first equation yields $x_4=\frac{a_2a_3x_1^2x_2}{a_1}$. Substituting 
this into the second equation yields $(x_1x_2x_3)^2=(\frac{a_1}{a_2a_3})^2\frac{b_1}{b_2b_3}$. 

%\textcolor{red}{Here assume $q$ odd!} Call $S$ the set of squares in $k^\times$. 
%Assume that 
If $\frac{b_1}{b_2b_3} \notin \FF_{q,2}^\times$, then the equalities are not both satisfied at the same time, hence 
the sum is zero. Let us then assume $\frac{b_1}{b_2b_3} \notin \FF_{q,2}^\times$. %, and let $\beta$ one of its square roots. 
Let us put $A:=\frac{a_1}{a_2a_3}$ and $B:=\frac{b_1}{b_2b_3}$. Then we have that 
$$(x_1, x_2, x_3, x_4) \in 
%T, \qquad T := 
\left\{(t, u, \frac{A\zeta}{tu}, \frac{t^2u}{A}) \mid t, u \in k^\times, \zeta^2=B\right\}.$$
The above sum can
%, which we call $S_{0,0}^{0,0}$ in Table \ref{tab:B2}, can
now be written as %split in two similar summands, one of which is 
\begin{align*}
&\sum_{\zeta^2=B,t, u \in k^\times} %\sum_{t, u \in k^\times} 
\phi(-t-\frac{A\zeta}{tu}-\frac{1}{a_3t}-\frac{tu}{a_2A\zeta}-\frac{t^3u}{a_2A^2\zeta}-\frac{t}{a_2A} 
-u-\frac{t^2u}{A}-\frac{t^2u}{b_3AB}-\frac{t^4u}{b_3A^2B}-2\frac{t^2}{b_3A\zeta}-\frac{1}{b_3u}-\frac{A}{b_2t^2u}),
\end{align*}
which can be put in the following more compact form, 
\begin{align*}
\sum_{t \in k^\times} \phi\left(-\frac{2t^2}{b_3A\zeta}-t\frac{a_2A+1}{a_2A}-\frac{1}{a_3t}\right)
%\cdot&[
%\sum_{u \in k^\times}
\calS_{q-1}\left(1,
%\phi(
%u(
-\frac{t^4}{b_3A^2\zeta^2}-\frac{t^3}{a_2A^2\zeta}-t^2\frac{b_3\zeta^2+1}{b_3A\zeta^2}-\frac{t}{a_2A\zeta}-1,%),
%u^{-1}(
-\frac{1}{b_3}-\frac{A\zeta}{t}-\frac{A}{b_2t^2}
\right)
%)
%]
.
\end{align*}
This sum, which we cannot reduce further, is the most complicated sum appearing in Table \ref{tab:B2}. All other sums in Tables \ref{tab:A2} and \ref{tab:B2} %e could not express directly 
can be expressed in a compact form, namely by at most four summands which only involve $\phi$ and the generalized Kloosterman sums 
$\calS_{q-1}(B, a, b)$ and $\tilde{\calS}_{q-1}(B, a, b, a', b')$.

\section{Structure constants in types $A_2$ and $B_2$}
\label{constants}
In this section we summarize our findings. In their action on the left, the elements of the standard basis of $\calH$ are given in the sequel. These elements are $e_0=e_0(a,b)$, $e_1=e_1(c)$, $e_2=e_2(d)$ and $e_3=1$ (as in \S\ref{basis_rank_2}). It is clear that $[e_3]=[e_3e_j:e_i]_{i,j=\overline{0,3}}$ is just the identity matrix.

We denote by $S_{ij}^{k}=[e_ie_j:e_k]$ the structure constants. Clearly $S_{ij}^{k}$ will depend on certain parameters: we use $a_1,b_1,c_1,d_1\in k^{\ast}$ for the index $i$, we use $a_2$, $b_2$, $\dots$ for the index $j$ and $a_3$, $b_3$, $\dots$ for the index $k$. For example $S_{02}^{1}=S_{02}^{1}(a_1,b_1,c_3,d_2)$.
Since $\calH$ is abelian, $S_{ij}^{k}=S_{ji}^{k}$ and we mark this redundant information with $\ast$ if $i<j$. %For example  $S_{21}^{0}$ in the matrix $[e_2]$.

The structure constants in types $A_2$ and $B_2$ are given in the second column of Tables \ref{tab:A2} and \ref{tab:B2} respectively. The third columns contain the distinguished subexpressions which parametrize the decompositions of $Un_iU\cap n_kUn_j^{-1}U$ as disjoint unions of left $U$-cosets (as in \S\ref{decomposition}). The last column gives the types of the distinguished subexpressions (as defined in \S\ref{sub:pt1}). 
We recall that all details for the computations in Tables \ref{tab:A2} and \ref{tab:B2}, most of which have been automated via the algorithm in Section \ref{algo}, can be found in \cite{PS}. 

%The easier to obtain structure constants, 
The structure constants which are easier to obtain, namely those of the form $q^{\ell(w)}$, are %included in the matrices of the elements. These are 
the ones corresponding to the $S_{ij}^{3}$ with $w_i^{-1}=w_j$. For such elements $w_i$ there is only one distinguished subexpression, namely the entire fixed reduced expression of $w_i$. Moreover, all its entries are of type A. For example, consider $S_{21}^{3}$ which is a sum over left $U$-coset representatives in
$
Un_2U\cap n_3Un_1^{-1}U
$.
The decomposition as in \S\ref{decomposition} has only one element, namely $J(w_2,w_1,w_3)=\{[2,1]\}$ corresponding to the fixed expression $s_2s_1$ of $w_2$. The type of this expression is AA. Moreover, it is easy to see from \eqref{eq:GGR} that all the summands are $\psi(1)=\phi(0)=1$.

We also recall that the $0$'s in the subexpressions mean that the corresponding simple reflection is omitted, e.g. $[1,2,0]=s_1s_2$ and $[1,0,1]=s_1s_1$ are subexpressions of $[1,2,1]=s_1s_2s_1$.

%\newpage

\subsection{Type $A_2$} In type $A_2$ there is an obvious extra symmetry to the structure constants obtained by interchanging the two simple roots, i.e. interchanging the indices $1$ and $2$ in $S_{ij}^{k}$.
$$
%[c_0(a)]_{i,j=\overline{0,3}}=
[e_0]
=
\begin{bmatrix}
S_{00}^{0} & S_{01}^{0} & S_{02}^{0} & \delta_{a_1,a_3}\delta_{b_1,b_3}%\delta_{(b_1, a_1), (b_3, a_3)}
\\
S_{00}^{1} & S_{01}^{1} & S_{02}^{1} & 0
\\
S_{00}^{2} & S_{01}^{2} & S_{02}^{2} & 0
\\
%q^3\delta_{(a_1, a_2), (b_2, b_1)} & 0 & 0 & 0
q^3\delta_{a_1, b_2}\delta_{a_2, b_1} & 0 & 0 & 0
%q^{\ell(w_0)} & 0 & 0 & 0
%q^{3} & 0 & 0 & 0
\end{bmatrix}
,%\quad
$$
$$
[e_1]=
\begin{bmatrix}
* & S_{11}^{0} &            S_{12}^{0} & 0
\\
* &           0 &            0 & \delta_{c_1, c_3}
\\
* & S_{11}^{2}  &           0 & 0
\\
0 &           0 &            q^2\delta_{c_1, -d_2} 
& 0
\end{bmatrix}
%$$
%Since $H$ is abelian, some entries were already obtained for $c_0(a,b)$ so the entries $*$ should be calculated only to check calculations.
%$$
%    [c_3]_{i,j=\overline{0,3}}=
,\quad
[e_2]=
\begin{bmatrix}
* &           * &          S_{22}^{0} & 0
\\
* &           0 &          S_{22}^{1} & 0
\\
* &           0 &          0 & \delta_{d_1, d_3}
\\
0 &  q^2\delta_{c_2, -d_1} 
&          0 & 0
\end{bmatrix}.
$$
\begin{table}[!ht]\label{tab:A2}
\renewcommand{\arraystretch}{1.15}
\begin{tabular}{ c c c c}
  \hline
  $S_{00}^{0}$
  &
  $
  \sum_{\zeta \mid \zeta^3=\frac{a_1b_1^2}{a_2^2a_3b_2b_3^2}}
\calS_{q-1}
\left(1,
\sigma_1,
\sigma_2
\right)$
 & [0,0,0] & BBB
 \\
 &
 $
  +q\delta_{a_1a_2b_3, -a_3b_1b_2}\delta_{a_1b_3, b_2b_3-b_1b_2}
 %+q\delta_{b_1b_2,a_1b_2-a_1b_3}
 $
&
[1,0,1]
&
CBA
  \\
  \hdashline
  $S_{00}^{1}$
  &
  $q\mathcal{S}_3\left(
\frac{a_1^2b_1}{a_2b_2^2}\frac{1}{c_3}, 
-1-\frac{b_2}{a_1}-\frac{a_2b_2}{a_1b_1},
-\frac{1}{b_2}-\frac{a_1}{a_2b_2}
\right)$
&
[0,2,0]
&
BAB
  \\
  $S_{00}^{2}$
  &
  $q\calS_3\left(\frac{a_1b_1^2}{a_2^2b_2d_3},-1-\frac{a_2}{b_1}-\frac{a_2b_2}{a_1b_1},-\frac{1}{a_2}-\frac{b_1}{a_2b_2}\right)$
  &
  [1,0,0]
&
ABB
  \\
  $S_{01}^{0}$
  &
  $
  \mathcal{S}_3
\left(
\frac{b_1c_2}{a_1a_3^2b_3}, 
1
+\frac{a_3}{b_1}, 
\frac{1}{a_1}
+\frac{1}{a_3}
+\frac{b_1}{a_3b_3}
\right)
  $
  &
  [0,2,0]
  &
BCB
  \\
  $S_{01}^{1}$
  &
  $q\delta_{a_1c_3,-c_2b_1}\phi\left(\frac{a_1}{c_2}\right)$
  &
  [1,0,1]
&
CBA
  \\
  $S_{01}^{2}$
  &
  $q\delta_{c_2d_3,-a_1b_1^2}\phi\left(\frac{b_1}{c_2}-\frac{b_1}{d_3}\right)$
    &
  [1,2,0]
&
ACB
  \\
  $S_{02}^{0}$
  &
  $
  \mathcal{S}_3\left(
\frac{a_3b_1}{a_1b_3d_2}, 
%1+\frac{b_3}{a_1}, \frac{1}{b_1}+\frac{1}{b_3}+\frac{a_1}{a_3b_3}
1+\frac{b_3}{a_1}, \frac{1}{b_1}+\frac{1}{b_3}+\frac{a_1}{a_3b_3}
\right)
  $
    &
  [0,0,1]
  &
BBC
  \\
  $S_{02}^{1}$
  &
  $q\delta_{c_3d_2,-a_1^2b_1}\phi\left(\frac{a_1}{d_2}-\frac{a_1}{c_2}\right)$
    &
  [0,2,1]
&
BAC
  \\
  $S_{02}^{2}$
  &
  $q\delta_{b_1d_3,-a_1d_2}\phi\left(\frac{b_1}{d_2}\right)$
    &
  [1,0,1]
&
ABC
  \\
  \hline
  $S_{11}^{0}$
  &
  $\delta_{c_1c_2,a_3^2b_3}\phi\left(\frac{a_3}{c_1}+\frac{a_3}{c_2}\right)$
    &
  [0,2]
&
BC
  \\
  $S_{11}^{2}$
  &
  $q\delta_{c_1,-d_3}\delta_{c_1,c_2}$
    &
  [1,2]
&
AC
  \\
  $S_{12}^{0}$
  &
  $\delta_{a_3c_1,b_3d_2}\phi\left(\frac{a_3}{d_2}\right)$
    &
  [1,0]
&
CB
  \\
  \hline
  $S_{22}^{0}$
  &
  $\delta_{d_1d_2,a_3b_3^2}\phi\left(\frac{b_3}{d_2}+\frac{b_3}{d_1}\right)$
    &
  [0,1]
&
BC
  \\
  $S_{22}^{1}$
  &
  $q\delta_{d_1,-c_3}\delta_{d_1,d_2}$
    &
  [2,1]
&
AC\\
\hline
\end{tabular}
\caption{Structure constants of $\calH$ when $G=\PGL_3(q)$}
\end{table}
%\end{TAB}
%\FloatBarrier

In the row corresponding to $S_{00}^0$ in Table \ref{tab:A2}, we put 
$$
\sigma_1=
-1
-\frac{b_3}{b_1}
+\zeta\left(
-\frac{a_2b_2b_3}{a_1b_1}%\zeta
-\frac{a_2b_3}{b_1}%\zeta
\right)
+\frac{1}{\zeta}\left(
-\frac{1}{a_2}
-\frac{1}{a_3}
\right), 
%-\frac{1}{a_2\zeta}
%-\frac{1}{a_3\zeta}, 
\qquad 
\sigma_2=
-\frac{1}{b_3}
-\zeta
-\frac{b_1}{a_2b_2b_3}\frac{1}{\zeta}.
$$

\subsection{Generation in type $A_2$}
\label{sub:genA}

We prove that $e_3, e_2(d)$ and 
$e_1(c)$ generate $\mathrm{A}_2$. By looking at the constant $S_{1 2}^0$ as in our table, we get 
\begin{align*}
e_1(c_1)e_2(d_2)-
q^2\delta_{c_1, -d_2}e_3&=\sum_{a_3, b_3 \in \FF_q^\times} 
\delta_{a_3c_1, b_3d_2}\phi(a_3/b_2)
e_0(a_3, b_3)
=\sum_{t \in \FF_q^\times} \phi(t)
e_0(d_2t, c_1t)
.
\end{align*}

In order to obtain just one $e_0(c, d)$ 
on the right hand side, we exploit the properties (i)--(iii) of $\phi$ in \S\ref{sub:firG}. We evaluate the above 
at $c_1z^{-1}$ and $d_2z^{-1}$ respectively, and we obtain
$$
e_1(c_1z^{-1})e_2(d_2z^{-1})-
q^2\delta_{c_1z^{-1}, -d_2z^{-1}}e_3
=
\sum_{t \in \FF_q^\times} \phi(t)
e_0(d_2tz^{-1}, c_1tz^{-1}).
$$
Notice that $\delta_{c_1z^{-1}, -d_2z^{-1}}=\delta_{c_1, -d_2}$. %We now look for a suitable function $f(z)$. 

Let us consider a function $f:\FF_q \to \CC$. By multiplying the above equation by $f(z)$ 
and then summing over $z$, we obtain
\begin{align*}
\sum_z f(z)\left(
e_1(c_1z^{-1})e_2(d_2z^{-1})-
q^2\delta_{c_1, -d_2}e_3
\right)
=
\sum_{t}\sum_{z} f(z)\phi(t)
e_0(d_2tz^{-1}, c_1tz^{-1}).
\end{align*}
We can now bring the term in $e_0$ outside 
one sum by the change of variable $s=t/z$, namely 
%the above can be written as 
$$
\sum_{s}\sum_{z} f(z)\phi(sz)
e_0(d_2s, c_1s)=\sum_{s}e_0(d_2s, c_1s)\sum_{z} f(z)\phi(sz).
$$
The claim is proved if we can choose the function $f$ %$f: \FF_q \to \CC^\times$ 
such that it satisfies the following property, 
$$\sum_{z} f(z)\phi(sz)=c_s\delta_{s, -1}
$$
for each $s \in \FF_q^\times$ and $c_s \in \FF_q^\times$; namely in this case 
the above expression can be written as 
$$\sum_{s}e_0(d_2s, c_1s)\sum_{z} f(z)\phi(sz)
=\sum_{s}e_0(d_2s, c_1s)c_s\delta_{s, -1}
=-e_0(d_2, c_1).
$$
%where $c=c_{-1}$. 

A function $f$ satisfying the above property does exist. Namely 
if $f(z)=\phi(z)-1$ then
$$\sum_{z} (\phi(z)-1)\phi(sz)=
\sum_{z}\phi((s+1)z)-\sum_{z}\phi(sz)
=(q-1)\delta_{s,-1}-\delta_{s \ne -1}+1=
q\delta_{s,-1}.
$$
%The claim is now proved. 
In particular we have obtained more explicitly, for each $x, y \in \FF_q^\times$,
%In particular we obtain, more explicitly, 
$$
e_0(x, y)=
q^{-1}
\sum_{z \in \FF_q^\times} 
(\phi(z)-1)e_1(-yz^{-1})e_2(-xz^{-1})
+
q^2\delta_{x, -y}e_3.
$$

\subsection{Type $B_2$} We assume here that $p\neq 2$.
$$
%[c_0(a)]_{i,j=\overline{0,3}}=
[e_0]
=
\begin{bmatrix}
%\sum_{i=1}^{3}S_{00}^{0,i}
S_{00}^{0}
& S_{01}^{0} &
%S_{02}^{0,1}+S_{02}^{0,2}
S_{02}^{0}
& \delta_{a_1, a_3}\delta_{b_1, b_3}
\\
%S_{00}^{1,0}+S_{00}^{1,1}
S_{00}^{1}
& S_{01}^{1} & S_{03}^{1} & 0
\\
S_{00}^{2} & S_{01}^{2} & S_{02}^{2} & 0
\\
q^4\delta_{a_1, a_2}\delta_{b_1, b_2} 
& 0 & 0 & 0
\end{bmatrix}
,
$$
$$
[e_1]=
\begin{bmatrix}
* & S_{11}^{0} &            S_{12}^{0} & 0
\\
* & 0            &              S_{12}^{1} & \delta_{c_1, c_3}
\\
* & S_{11}^{2}  &          S_{12}^{2} & 0
\\
0 & q^3\delta_{c_1, c_2} 
&           0 & 0
\end{bmatrix}
%$$
%Since $H$ is abelian, some entries were already obtained for $c_0(a,b)$ so the entries $*$ should be calculated only to check calculations.
%$$
%    [c_3]_{i,j=\overline{0,3}}=
,\quad
[e_2]=
\begin{bmatrix}
* &  * &           S_{22}^{0} & 0
\\
* & * &           S_{22}^{1} & 0
\\
* & * &           0   & \delta_{d_1, d_3}
\\
0 & 0 &           q^3\delta_{d_1, -d_2} 
& 0
\end{bmatrix}
$$
%\FloatBarrier
\begin{table}[!h]\label{tab:B2}
\renewcommand{\arraystretch}{1.2}
\begin{tabular}{ c c c c c}
  \hline
  $S_{00}^{0}$
  & determined in \S\ref{sub:S000}
  &
  \\
  \hdashline
  $S_{00}^{1}$
  &
  $
  q\sum_{\zeta^2=\frac{b_1}{b_2c_3}}\calS_{q-1}(1, 
  -1
-\frac{a_1}{a_2}
-\frac{b_2}{b_1}
-\frac{a_2b_2}{a_1b_1}
-\frac{1}{a_2\zeta}
-\frac{1}{a_1\zeta}, 
\zeta
-\frac{1}{b_2})
$
   &
   [1,0,0,0]
   &
   ABBB
\\
&
$q^2\delta_{a_1, -a_2}\sum_{\zeta^2=\frac{b_1}{b_2c_3}}\delta_{\frac{1}{a_2\zeta}, 1-\frac{b_1}{b_2}}$
%$
   %+ q^2\delta_{a_1, -a_2}(\delta_{c_3, \frac{(b_1-b_2)^2}{b_1b_2}a_2^2}+\delta_{c_3, -\frac{(b_1-b_2)^2}{b_1b_2}a_2^2})
   %$
      &
   [1,2,0,2]
   &
   ACBA
  \\
  \hdashline
  $S_{00}^{2}$
  &
  $
  q\sum_{\zeta^2=\frac{b_1}{b_2}}\tilde{S}_{q-1}(1, \tau_1', \tau_2', \tau_3', \tau_4')
  $
  &
  [0,2,0,0]
  &
  BABB
  \\
  &
  $
  +q\delta_{b_1,b_2}(2\delta_{-\frac{a_2d_3}{a_1} \in \FF_{q,2}}-1)G\phi(\frac{d_3}{4a_1a_2})
  $
  &
  [1,0,1,2]
  &
  CBAA
  \\
  \hdashline
  $S_{01}^{0}$
  &
  $
  \sum_{\zeta^2=\frac{c_2}{b_1b_3}}\calS_{q-1}(1, -\zeta b_3-\frac{a_3b_3}{a_1}\zeta -\frac{c_2}{b_1}-\frac{a_3c_2}{a_1b_1},
  -\frac{b_1}{a_3c_2}\zeta-\frac{a_1b_1}{a_3b_3c_2}-\frac{1}{c_2}
  )
  $
  &
  [0,0,1,0]
  &
  BBCB
  \\
  &
  $
  +
  q\delta_{a_1, -a_3}\sum_{\zeta^2=\frac{b_1}{b_3c_2}}\delta_{\frac{1}{\zeta}, -\frac{c_2}{b_1}}
  %q\delta_{a_1, -a_3}(\delta_{c_2, b_1b_3}+\delta_{c_2, -b_1b_3})
  $
  &
  [1,2,0,2]
  &
  CCBA
  \\
  \hdashline
  $S_{01}^{1}$
  &
  $q\calS_2(b_2c_2c_3, 0, a_1b_1+c_2+c_3)$
  %$q\sum_{\zeta}\phi(\frac{c_2+a_1b_1+c_3}{\zeta})$
  %&
  %$\zeta^2=c_2c_3b_2$
  &
  [1,0,1,0]
  &
  ABCB
  \\
  \hdashline
  $S_{01}^{2}$
  &
  $q\phi(-\frac{a_1}{d_3}-\frac{d_3}{a_1b_1})\calS_2(\frac{b_1}{c_2}, 1, \frac{1}{a_1})$
  %$q\sum_{\zeta}\phi(\zeta + \frac{1}{a_{1}\zeta}-\frac{d_3}{a_1b_1} - \frac{a_1}{d_3})$
  &
  %$\zeta^2=\frac{b_1}{c_2}$
  %&
  [0,2,1,0]
  &
  BACB
  \\
  \hdashline
  $S_{02}^{0}$
  &
  $\sum_{\zeta^2=\frac{b_1}{b_3}}\tilde{\calS}_{q-1}(1, \tau_1, \tau_2, \tau_3, \tau_4)$
  &
  [0,0,0,2]
  &
  BBBC
  \\
  &
  $+\delta_{b_1,b_3}(2\delta_{\frac{a_3d_2}{a_1b_3} \in \FF_{q,2}}-1)G\phi(-\frac{a_1b_3}{4a_3d_2}(1-\frac{a_3}{a_1})^2)$
  &
  [1,0,1,2]
  &
  CBAC
  \\
  \hdashline
  $S_{02}^{1}$
  &
  $q\phi(\frac{a_1}{d_2}+\frac{d_2}{a_1b_1})\calS_2(\frac{b_1}{c_3}, 1, \frac{1}{a_1})$
  %$q\sum_{\zeta}\phi(\zeta+\frac{1}{a_1\zeta}+\frac{a_1}{d_2} + \frac{d_2}{b_1a_1})$
  %&
  %$\zeta^2=\frac{b_1}{c_3}$
  &
  [1,0,0,2]
  &
  ABBC
  \\
  \hdashline
  $S_{02}^{2}$
  &
  $q\sum_{\zeta_1^2=-\frac{a_1}{d_2d_3}}\sum_{\zeta_3^2=-\frac{a_1b_1}{d_2d_3}}\phi\left(
\zeta_3-\frac{1}{d_3\zeta_1}+\frac{1}{d_2\zeta_1}+2\frac{\zeta_1}{\zeta_3}
\right)$
  %%%$q\sum_{\zeta_1^2=b_1,\zeta_2^2=-\frac{a_1}{d_2d_3}}\phi(\zeta_1 \zeta_2 - \frac{1}{d_{3}\zeta_2}
  %+ \frac{1}{d_{2}\zeta_2}+2\frac{1}{\zeta_1})$%%%
  %&
  %$\zeta_1^2=b_1,\zeta_2^2=-\frac{a_1}{d_2d_3}$
  &
  [0,2,0,2]
  &
  BABC
  \\
    \hline
  $S_{11}^{0}$
  &
  $\mathcal{S}_2(\frac{c_1}{b_3c_2}, 1, \frac{a_3}{c_2}+\frac{1}{b_3})$
  %$\calS_2(b_3c_1c_2, 0, a_3b_3+c_1+c_2)$
  %$\sum_{\zeta}\phi\left(\frac{c_1+a_3b_3+c_2}{\zeta}\right)$
  %&
  %$\zeta^2=b_3c_1c_2$
  &
  [0,1,0]
  &
  BCB
  \\
  $S_{11}^{2}$
  &
  $q\delta_{c_1, c_2}\phi\left(-\frac{d_3}{c_2}\right)$
  %&
  %$c_1=c_2$
  &
  [2,1,0]
  &
  ACB
  \\
  $S_{12}^{0}$
  &
  $\phi\left(\frac{a_3}{d_2} + \frac{d_2}{b_{3}a_{3}} \right)\calS_2(\frac{b_3}{c_1}, 1, \frac{1}{a_3})$
  %$\phi\left(- \frac{1}{a_{3}\zeta}\frac{a_3}{d_2} + \frac{d_2}{b_{3}a_{3}}\right)$
  %$\phi(\zeta + \frac{1}{a_{3}\zeta}+\frac{a_3}{d_2} + \frac{d_2}{b_{3}a_{3}})+\phi(-\zeta - \frac{1}{a_{3}\zeta}\frac{a_3}{d_2} + \frac{d_2}{b_{3}a_{3}})$
 %&
 %$\zeta^2=\frac{b_3}{c_1}$
  &
  [0,0,2]
  &
  BBC
  \\
  $S_{12}^{1}$
  &
  $q\delta_{c_1,c_3}\phi\left(\frac{d_2}{c_3}\right)$
  %&
  %$c_1=c_3$
  %&
  &
  [0,1,2]
  &
  BAC
  \\
  $S_{12}^{2}$
  &
  $q\delta_{d_2,-d_3}\sum_{\zeta^2=c_1} \phi\left(\frac{\zeta}{d_3}\right)$
  %$\delta_{d_2,-d_3}\calS_2(c_1,\frac{1}{d_3},1)$
  %$\phi\left(\frac{\zeta}{d_3}\right)+\phi\left(-\frac{\zeta}{d_3}\right)$
%  &
%  $d_2=-d_3$ and $\zeta^2=c_1$
%  &
%  $q$
&
[2,0,2]
&
ABC
  \\
  \hline
  $S_{22}^{0}$
  &
  $\sum_{\zeta_1^2=\frac{d_1}{a_3d_2}}
\sum_{\zeta_3^2=\frac{d_1}{a_3b_3d_2}}
\phi\left(
-\zeta_1-\frac{1}{d_2\zeta_3}-\frac{1}{a_3\zeta_1}+2\frac{\zeta_1}{\zeta_3}
\right)$
  %$
  %q\sum_{\zeta_1^2=\frac{1}{b_3},\zeta_2^2=\frac{d_1}{a_3d_2}}
  %\phi\left(-\zeta_2 - \frac{1}{d_{2}\zeta_1\zeta_2} - \frac{1}{a_{3}\zeta_2}+2\zeta_1\right)$
  
  %&
 &
 [0,2,0]
 &
 BCB
 \\
   $S_{22}^{1}$
  &
  %$q\phi(\zeta)+q\phi(-\zeta)$
  $q\delta_{d_1, d_2}\calS_2\left(\frac{1}{c_3},0, \frac{1}{d_2}\right)$
  &
  [1,2,0]
  &
  ACB
%  &
%  $\zeta=\frac{d_1d_2}{c_3}$
%  &
%  $2q$
\\
\hline
\end{tabular}
\caption{Structure constants of $\calH$ when $G=\SO_5(q)$ and $p$ is odd}
\end{table}
%\FloatBarrier

In the rows corresponding to $S_{02}^{0}$ and $S_{00}^{2}$ in Table \ref{tab:B2}, we put 
$$
\tau_1=a_3d_2(\frac{2}{a_1b_3\zeta}+\frac{1}{a_1b_3}+\frac{1}{a_1b_1}), 
 \qquad
\tau_2=-1-\zeta -\frac{a_3}{a_1\zeta}-\frac{a_3}{a_1},
 \qquad
\tau_3=-\frac{1}{a_3}, \qquad \tau_4=\frac{a_1}{a_3d_2},
$$
%and
$$
\tau_1'=a_2d_3(\frac{2}{a_1\zeta}-\frac{b_2}{a_1b_1}-\frac{1}{a_1}),
\qquad
\tau_2'=d_3(\frac{1}{a_1\zeta}-\frac{1}{a_1}), 
\qquad 
\tau_3'=\frac{a_1\zeta}{a_2d_3}-\frac{1}{d_3}, 
\qquad
\tau_4'=-\frac{a_1}{a_2b_2d_3}.
$$

\subsection{Generation in type $B_2$}
The problem of determining whether each of the $e_0(a, b)$ is generated by elements of the form $e_1(c)$, $e_2(d)$ and $e_3$ 
seems much harder in this case. Here we again have %In fact, by looking at $S_{12}^0$, we again have that 
%Notice the following in type $\mathrm{B}_2$. By looking at the constant $S_{1 2}^0$, we have that 
\begin{equation}\label{eq:12}
e_1(c_1)e_2(d_2)+K
=
\sum_{a_3, b_3}S_{12}^0e_0(a_3, b_3)%+K, \qquad  \in \langle\langle e_1(c), e_2(d), e_3 \mid c, d \in \FF_q \rangle\rangle.
\end{equation}
for some element $K$ generated by $e_1(c)$, $e_2(d)$ and $e_3$. %Applying the same methods as in 
Notice that if we let $\langle \mu \rangle=\FF_q^\times$ then we can write
%We look at the right hand side of \eqref{eq:12}. Let $\langle \mu \rangle=\FF_q^\times$. We have that 
\begin{align*}
\sum_{a_3, b_3}S_{1 2}^0e_0(a_3, b_3)
&=\sum_{t}\sum_{i=1}^{\frac{q-1}{2}}
\left(
\phi(t + \frac{1}{c_1\mu^{2i}t}
+\mu^i + \frac{1}{d_2\mu^it})+
\phi(t + \frac{1}{c_1\mu^{2i}t}
-\mu^i - \frac{1}{d_2\mu^it})
\right)
e_0(d_2t, c_1\mu^{2i})\\
&=\sum_{t}\sum_{u}
\left(
\phi(t + \frac{1}{c_1tu^2}
+u + \frac{1}{d_2tu})
\right)
e_0(d_2t, c_1u^2)
\end{align*}
since $\mu^{(q-1)/2}=-1$. 

%Now 
We try to apply the same method as in \S\ref{sub:genA}. Let us consider a function $f: \FF_q^2 \to \CC$. By replacing $c_1$ and $d_2$ in \eqref{eq:12} with $c_1z^{-2}$ and $d_2w^{-1}$ respectively, 
multiplying by $f(z, w)$, and summing over $z$ and $w$, we get 
\begin{align*}
\sum_{z, w}f(z, w)
e_1(c_1z^{-2})e_2(d_2w^{-1})+K'
&=\sum_{x, y}e_0(d_2x, c_1y^2)\sum_{z, w}f(z, w)
\phi
\left(
zy+wx+\frac{1}{d_2xyz}+\frac{1}{c_1wxy^2}
\right),
\end{align*}
The problem of finding a suitable choice of $f$ which allows to express %expressing 
the second sum on the right hand side as a product of suitable 
delta functions, which would leave just one term $e_0(x, y)$, remains an open question.

\end{document}